\documentclass[11pt,twoside]{article}

\usepackage{latexsym}
\usepackage{amsmath}
\usepackage{amssymb}
\usepackage{amscd}
\usepackage{eucal}
\usepackage{amsthm}

\newtheorem{thm}{Theorem}[section]
\newtheorem{prop}[thm]{Proposition}
\newtheorem{lemm}[thm]{Lemma}
\newtheorem{defi}[thm]{Definition}

\DeclareMathOperator{\Ext}{Ext}
\DeclareMathOperator{\projdim}{projdim}
\begin{document}
\title{Canonical surfaces in $\mathbb{P}^4$ and Gorenstein algebras in codimension 2}
\author{Christian B\"ohning (G\"ottingen)}
\date{}
\maketitle
\begin{quote}
\small
\begin{center}\bfseries CONTENTS\\ \end{center}
\vspace{0.5cm}
\mdseries
0. Introduction \dotfill\ p. 01\\
1. R\'esum\'e of properties of canonical surfaces in $\mathbb{P}^4$ \dotfill\ p. 04\\
2. Analysis of the case $K^2=11$ \dotfill\ p. 11\\
3. The case of general $K^2$ \dotfill\ p. 22\\
4. A commutative algebra lemma \dotfill\ p. 32\\
\end{quote}
\normalsize
\setcounter{section}{-1}
\section[Introduction]{Introduction}
There are 2 major circles of problems providing motivational background for this
work:\\ The first centers around the questions in the theory of algebraic surfaces of
general type concerning the existence of surfaces with prescribed invariants (the
"geography") and, more systematically, the problem of describing their moduli spaces
and canonical resp. pluricanonical models; I want to consider this general problem
in the special case of surfaces with geometric genus $p_g=5$, more precisely for
canonical surfaces in $\mathbb{P}^4$ (i.e. those for which the canonical map is a
birational morphism) with $q=0$ and $p_g=5$.\\
The second type of questions is more algebraic in spirit with principal aim to find a
satisfactory structure theorem for Gorenstein algebras in codimension 2; roughly,
these are finite $R-$algebras $B$ ($R$ some "nice" base ring) with
$B\cong\mathrm{Ext}^2_R(B,R)$ (cf. section 1 below for precise definitions). With
regard to a structure theorem, "satisfactory" means that one should be able to tell
from practically verifiable and non-tautological conditions how the Hilbert
resolution of $B$ over $R$ encodes 1) the "duality" $B\cong\mathrm{Ext}^2_R(B,R)$ and
2) the fact that $B$ has not only an $R-$module structure, but also a ring structure.
Whereas 1) is by now fairly well understood, 2) is not.\\
The link between these two problems is the following: Given a regular surface $S$ of
general type with canonical map $S\to Y\subset \mathbb{P}^4$ a birational morphism,
the canonical ring $\mathcal{R}=\bigoplus\limits_{m\ge 0}H^0(S,\mathcal{O}_S(mK))$ is
a codimension 2 Gorenstein algebra over $\mathcal{A}=\mathbb{C}[x_0,\ldots,x_4]$, the
homogeneous coordinate ring of $\mathbb{P}^4$. And conversely, starting from a
Gorenstein algebra $\mathcal{R}$ in codimension 2 over $\mathcal{A}$, one recovers in
shape of $X=\mathrm{Proj}(\mathcal{R})$ the canonical model of a surface $S$ as
above, provided
$\mathrm{Ann}_{\mathcal{A}}(\mathcal{R})$ is prime
(i.e. $Y=\mathrm{Supp}(\mathcal{R})\subset
\mathbb{P}^4$ with its reduced induced closed subscheme structure is an irreducible
surface), 
$X$ has only rational double points as singularities, and some weak technical
assumption on a presentation matrix of $\mathcal{R}$ as $\mathcal{A}-$module holds
(cf. thm. 1.6 below).\\

 Let me give a little history of the development that lead to
these ideas. The structure of canonical surfaces in
$\mathbb{P}^4$ with
$q=0,\:p_g=5,\:K^2=8$ and 9 was worked out already by F. Enriques (cf. [En], p.
284ff.; they are the complete intersections of type (2,4) and (3,3)). The case
$K^2=10$ was solved by C. Ciliberto using liaison arguments. However, these techniques
could not be utilized for higher
$K^2$. D. Ro\ss berg has constructed examples for the cases $K^2=11$ and $=12$, and
gives a partial description of the moduli spaces of these surfaces. He constructs
these surfaces as degeneracy loci of morphisms between reflexive sheaves of rank $n$
and $n+1$. The approach taken in this work relies instead on ideas in [Cat2]. In this
paper canonical surfaces in $\mathbb{P}^3$ are studied (from a moduli point of view)
via a structure theorem proved therein for Gorenstein algebras in codimension 1. It
is shown that the duality $B\cong \mathrm{Ext}^1_R(B,R)$ for these algebras
translates into the fact that the Hilbert resolution of $B$ can be chosen to be
self-dual; moreover, that the presence of a ring structure on $B$ is equivalent to a
(closed) condition on the Fitting ideals of a presentation matrix of $B$ as
$R-$module (the so-called "ring condition" or "rank condition" or "condition of
Rouch\'e-Capelli", abbreviated R.C. in any case) whence the moduli spaces of the
surfaces studied in that paper can be parametrized by locally closed sets of
matrices. These ideas were developed further and generalized in [M-P] and [dJ-vS]
(within the codimension 1 setting). In particular the latter papers show that R.C.
can be rephrased in terms of annihilators of elements of $B$ and gives a good
structure theorem also in the non-Gorenstein case. Subsequently, M. Grassi isolated in
[Gra] the abstract kernel of the problem and proved that also for codimension 2
Gorenstein algebras the duality $B=\mathrm{Ext}^2_R(B,R)$ is equivalent to $B$ having a
self-dual Hilbert resolution. He introduced the concept of Koszul modules which
provide a nice framework for dealing with Gorenstein algebras and also proposed a
structure theorem for the codimension 2 case. Unfortunately, as for the question of
how the ring structure of a codimension 2 Gorenstein algebra is encoded in its
Hilbert resolution, the conditions he gives are tautological and (therefore) too
complicated (although they are necessary and sufficient). More recently, D. Eisenbud
and B. Ulrich ([E-U]) re-examined the ring condition and gave a
generalization of it which appears to be more natural than the direction in which
[Gra] is pointing. But
 essentially, they only give sufficient conditions for $B$ to be a ring, and these
are not fulfilled in the applications to canonical surfaces one has in mind. More
information on the development sketched here can be found in [Cat4].  For a deeper
study of that part of the story that originates from the duality $B\cong
\mathrm{Ext}^2_R(B,R)$ and its effects on the symmetry properties of the Hilbert
resolution of $B$, as well as for a generalization of this to the bundle case cf.
[Wal].\\
The aim of this work is to show that whereas in codimension 1, R.C. and
the symmetry coming from the Gorenstein condition can be treated as separate
concepts, in codimension 2 they seem to be more intimately linked: For canonical
surfaces in $\mathbb{P}^4$ with $p_g=5,\: q=0, \: K^2\ge 10$ the symmetry implies
R.C. under some mild extra conditions. Moreover, the ideas in [E-U] can be made to work
in the geometric setting with some additional effort. Finally, a slightly stronger
statement than in [Gra] can be made concerning the structure of the Hilbert
resolution of a codimension 2 Gorenstein algebra. 
\\
The single sections are organized as follows:\\
Section 1 is meant to establish the connection mentioned above between canonical
surfaces in
$\mathbb{P}^4$ and Gorenstein algebras in codimension 2 and thus contains no new 
results. Some proofs have been added, partly for completeness, partly because
occasionally minor simplifications of known proofs could be made by pulling
together results from different sources, and sometimes adaptations to the chosen
set-up were necessary.\\
In section 2 I study how for canonical surfaces with $q=0,\: p_g=5,\: K^2=11$ the
ring structure of the canonical ring $\mathcal{R}$ is encoded in its Hilbert
resolution by reduction of a presentation matrix of $R$ to a normal form (modelled
on an example in [Ro\ss ], p. 115). I show how this can be used to derive
information on the moduli space of these surfaces.\\
In section 3 I use a localization argument to reduce, to a major extent, the study of
the case of general $K^2$ to structural aspects of the case $K^2=10$ which
is considered in [Cat4]. In particular it should be possible on the ground of this
to subsequently decide on the existence of these surfaces and to study their moduli
spaces (the problem largely boils down to an understanding of the Fano variety of
$\mathbb{P}^4$'s on certain varieties of "Gorenstein-symmetric complexes of length
2"). In particular, the automaticity of the ring condition for these surfaces is
established under the assumption that their images $Y\subset\mathbb{P}^4$ have
only improper double points in their nonnormal locus (a fact already hinted at in
[Cat4], p.  48). Whether this last condition is really needed, or only veils some 
algebraic counterpart behind it, is doubtful (cf. the remarks at the end of this
section).\\
Section 4 contains a result that improves on work in [Gra]; I have found no
further use of it by now, but it could be embodied as a technical lemma in the
attempt of finding a good structure theorem for Gorenstein algebras in
codimension 2.    
\\
As for standard notation from surface theory (such as $p_g$ for the geometric genus,
$P_m$ for the plurigenera etc.) I adhere to [Beau] (except that I write $e(S)$
instead of $\chi_{\mathrm{top}}(S)$ for the topological Euler characteristic of $S$).
My commutative algebra notation agrees largely with [Ei], but the following point
(which traditionally seems to cause notational confusion) should be noted: For
$I\subset R$ an ideal in a Noetherian ring and $M$ a finite $R-$module, I write
$\mathrm{grade}(I,M)$ for the length of a maximal $M-$regular sequence contained in
$I$ ($=\mathrm{min}\{ i:\:\mathrm{Ext}_R^i(R/I,M)\neq 0\}$), and also, if there is no
risk of confusion, $\mathrm{grade}\, M:=\mathrm{grade}(\mathrm{Ann}_R(M),R)$ and
$\mathrm{grade}\,I:=\mathrm{grade}(I,R)$. Furthermore if $R=(R,\mathfrak{m},k)$ is a
Noetherian local ring or graded ring with $\mathfrak{m}$ a unique maximal element
among the graded proper ideals of $R$ (e.g. a positively graded algebra
over a field), I write  $\mathrm{depth}\, R:=\mathrm{grade}(\mathfrak{m},R)$. This is
in accordance with [B-He] and the terminology seems to go back to Rees.

\section[R\'esum\'e of canonical surfaces]{R\'esum\'e of properties of canonical
surfaces in
$\mathbb{P}^4$} In this section I want to gather together the results on canonical
surfaces in
$\mathbb{P}^4$ needed in the sequel and give proofs for the more
important of them.
\begin{defi}
Let $S$ be a smooth surface and $\pi:S\to Y\subseteq \mathbb{P}^4$ a
morphism given by a 5-dimensional base-point free linear subspace $L$ of
$H^0(S,\mathcal{O}_S(K))$ and such that $\pi$ is birational onto its image
$Y$ in $\mathbb{P}^4$. Then Y is called a \underline{canonical surface} in
$\mathbb{P}^4$ (and $\pi$  an almost generic canonical projection).
\end{defi}
In the above situation, since $K_S$ is nef, $S$ is automatically a
minimal model of a surface of general type.\\ Henceforth I will make the
assumption that
$S$ is a regular surface, i.e.\\ \underline{
$q=h^1(S,\mathcal{O}_S)=0$}, basically because then the canonical ring
$\mathcal{R}:=$\\ $=\bigoplus\limits_{n\ge 0}H^0(S,\mathcal{O}_S(nK))$
enjoys the following property  which makes it convenient to study by
homological methods:
\begin{prop}
$\mathcal{R}$, viewed as a module over the homogeneous coordinate ring of
$\mathbb{P}^4\; \mathcal{A}=\mathbb{C}[x_0,\ldots,x_4]$ via $\pi$, is a
Cohen-Macaulay (CM) module iff S is a regular surface.\end{prop}
\begin{proof}Since $K_S$ is nef and big, the Ramanujam vanishing theorem
gives $H^1(S,\mathcal{O}_S(lK_S))=0$ for $l\in\mathbb{Z},\;l<0$, which
holds also for
$l\geq 2$  by Serre duality $H^1(S,\mathcal{O}_S(lK_S))\cong
H^1(S,\mathcal{O}_S((1-l)K_S))$ on $S$. Therefore S is regular iff
$H^1(S,\mathcal{O}_S(lK))=0\; \forall l\in \mathbb{Z}$ (taking again into
account $H^1(S,\mathcal{O}_S(K_S))\cong H^1(S,\mathcal{O}_S)$). I will
prove that the latter is equivalent to
$\mathcal{R}$ being CM.\newline In fact, if
$\mathcal{R}$ is CM, $\projdim_{\mathcal{A}}\mathcal{R}=2$ and in particular I have
that\\
$\Ext_{\mathcal{A}}^3(\mathcal{R}(i),\mathcal{A}(-5))=0\; \forall  i$. By Serre
duality
$H^1(\mathbb{P}^4,\widetilde{\mathcal{R}(i)})=0\;\forall i$, the tilde denoting
the sheaf associated to a graded module. But
$\tilde{\mathcal{R}}\cong \pi_{\ast}\mathcal{O}_S$, hence
$H^1(S,\mathcal{O}_S(lK_S))=0\;\forall l \in
\mathbb{Z}$.\\Conversely, suppose $H^1(S,\mathcal{O}_S(lK_S))=0\;\forall
l\in\mathbb{Z}$. The idea is now to derive the CM property of $\mathcal{R}$ by
looking at $C$ on $S$, the pullback of a generic hyperplane section $H$ of $Y$ via
$\pi$: genericity means that $C$ is a nonsingular divisor in $L$ such that $\pi |_C$ is
a birational morphism  onto $H$ (by Bertini's theorem such $C$ exists). It is known
that then $\mathcal{R}'=\bigoplus\limits_{m\ge 0}H^0(C,\mathcal{O}_C(mK_S))$ is CM (cf.
[Sern], lemma 1.1). Assume $H$ is cut out on $Y$ by $x_4=0$ and consider for each $m$ 
the cohomology long exact sequence
\begin{eqnarray*}
0 \longrightarrow H^0(S,\mathcal{O}_S((m-1)K_S)) \stackrel{\cdot
x_4}{\longrightarrow} H^0(S,\mathcal{O}_S(mK_S)) \longrightarrow
H^0(C,\mathcal{O}_C(mK_S)) \\ \longrightarrow H^1(S,\mathcal{O}_S((m-1)K_S))
\longrightarrow \ldots
\end{eqnarray*}
Since the $H^1-$terms vanish by hypothesis, I get $\mathcal{R}'=\mathcal{R}/
x_4 \mathcal{R}$, whence 
\begin{eqnarray*}
2=\dim (\mathcal{R}')=\mathrm{depth}_{\mathbb{C}[x_0,\ldots ,x_3]}(\mathcal{R}')=
\mathrm{depth}_{\mathbb{C}[x_0,\ldots ,x_3]}(\mathcal{R}/ x_4 \mathcal{R})\\ \le
\mathrm{depth}_{\mathbb{C}[x_0,\ldots ,x_4]}(\mathcal{R})-1
\end{eqnarray*}
since $x_4$ is regular on $\mathcal{R}$. But as $\dim (\mathcal{R})=3$ and
$\mathrm{depth}_{\mathbb{C}[x_0,\ldots ,x_4]}(\mathcal{R})\le \dim (\mathcal{R})$, I
get $\mathrm{depth}_{\mathbb{C}[x_0,\ldots ,x_4]}(\mathcal{R})= \dim (\mathcal{R})$,
i.e. $\mathcal{R}$ is CM.
\end{proof}
On the other hand, the fact that $\pi$ is an almost generic canonical
projection
($\mathcal{O}_S(K)\cong\pi^{\ast}\mathcal{O}_{\mathbb{P}^4}(1)$) implies
that various remarkable duality statements hold for $\mathcal{R}$, which
I shall frequently exploit in the sequel and which can be best expressed
in terms of properties of the minimal free resolution of $\mathcal{R}$.
Precisely:
\begin{defi}Let $R:=k[x_1,\ldots,x_r]$ be a polynomial ring in $r$
indeterminates over some field $k$, graded in the usual way, and let $B$
be a graded $R$-algebra. $B$ is said to be a \underline{Gorenstein algebra
of codimension $c$} (and with twist $d\in \mathbb{Z}$) over $R$ $:\Longleftrightarrow
B\cong\Ext^c_R(B,R(d))$ as B-modules.
\end{defi}
[The $B$-module structure on $\Ext_R^c(B,R(d))$ is induced from $B$
by functoriality of $\Ext_R^c(\cdot,R(d))$: If $b\in B$ and $m_b:B\to B$
is multiplication by $b$ on $B$, the map $\Ext_R^c(m_b,R(d))$ is
multiplication by $b$ on $\Ext_R^c(B,R(d))$.]
\begin{thm}
With the hypotheses and notation of definition 1.1 $\mathcal{R}$ is a
Gorenstein algebra of codimension 2 over $\mathcal{A}$ and as such has a
minimal graded free  resolution of the form:
\begin{eqnarray*} \begin{CD}\mathbf{R}_{\bullet}:\;
0 @>>> \bigoplus \limits_{i=1}^{n+1}\mathcal{A}(-6+r_i)
@>{{-\beta^t}\choose{\alpha^t}}>>
\bigoplus
\limits_{j=1}^{n+1} \mathcal{A}(-6+s_j)\oplus\bigoplus
\limits_{j=1}^{n+1}\mathcal{A}(-s_j)\end{CD} \\ \begin{CD}@>{(\alpha\,\beta)}>>
\bigoplus
\limits_{i=1}^{n+1}\mathcal{A}(-r_i) @>>> \mathcal{R} @>>> 0 
\end{CD}\, .
\end{eqnarray*}

\end{thm}
\begin{proof}\emph{(sketch)}
 Setting $X:=\mathrm{Proj}(\mathcal{R})$, the canonical model of $S$, I get that
$\pi$, being given by a base-point free linear subsystem of $|K_S|$, factors
through  $X$ as in the picture:\\
\setlength{\unitlength}{1cm}
\begin{center}
\begin{picture}(4,3)
\put(0.7,2.5){$S$}
\put(2.8,2.5){$Y\subset \mathbb{P}^4$}
\put(1.75,0.7){$X$}
\put(1,2.6){\vector(1,0){1.6}}
\put(0.8,2.3){\vector(2,-3){0.8}}
\put(2,1.1){\vector(2,3){0.8}}
\put(1.75,2.7){$\pi$}
\put(2.5,1.5){$\psi$}
\put(1,1.5){$\kappa$} 
\end{picture}\end{center}
and $\psi$ is a finite morphism onto $Y$. Hence by relative duality for finite
morphisms (cf. e.g. [Lip], p. 48ff.),
$\psi_{\ast}\omega_X=\mathcal{H}om_{\mathcal{O}_Y}(\psi_{\ast}\mathcal{O}_X,\omega_Y)$,
where
$\omega_Y=\mathcal{E}xt^2_{\mathcal{O}_{\mathbb{P}^4}}(\mathcal{O}_Y,\omega_{\mathbb{P}^4})$
 and $\omega_X$ are the Grothendieck dualizing sheaves of $Y,\: X$ resp.; but 
$\mathcal{H}om_{\mathcal{O}_Y}(\psi_{\ast}\mathcal{O}_X,\omega_Y)=\mathcal{E}xt^2_{\mathcal
{O}_{\mathbb{P}^4}}(\psi_{\ast}\mathcal{O}_X,\omega_{\mathbb{P}^4})$ since $Y$ has
codimension 2 in $\mathbb{P}^4$ (cf. also [Har], p. 242). Furthermore
$\psi_{\ast}\omega_X=\tilde{\mathcal{R}}(1)$ (cf. [Cat 2], p.76, prop. 2.7) and
$\psi_{\ast}\mathcal{O}_X=\tilde{\mathcal{R}}$. Thus I get
\begin{equation}
\tilde{\mathcal{R}}=\mathcal{E}xt^2_{\mathcal{O}_{\mathbb{P}^4}}(\tilde{\mathcal{R}},
\mathcal{O}_{\mathbb{P}^4}(-6)).
\end{equation}
Since $\mathcal{R}$ is CM I get a length 2 resolution
\begin{equation}
0\to F_2 \to F_1 \to F_0 \to \mathcal{R} \to 0,
\end{equation}
with $F_0,\: ,F_1,\: F_2$ graded free $\mathcal{A}-$modules.
Taking $\mathrm{Hom}_{\mathcal{A}}(\cdot ,\mathcal{A}(-6))$ of (2) I obtain
\begin{equation}
0 \to F^{\vee}_0(-6) \to F_1^{\vee}(-6) \to F^{\vee}_2(-6) \to
\mathrm{Ext}^2_{\mathcal{A}}(\mathcal{R},\mathcal{A}(-6)) \to 0.
\end{equation}
From the facts that $\mathcal{E}xt^2_{\mathcal{O}_{\mathbb{P}^4}}(\tilde{\mathcal{R}},
\mathcal{O}_{\mathbb{P}^4}(-6))$ is the sheaf associated to\\
$\mathrm{Ext}^2_{\mathcal{A}}(\mathcal{R},\mathcal{A}(-6))$, and
$\tilde{\mathcal{R}}$ the sheaf associated to $\mathcal{R}$, and I have resolutions
(2) and (3) of length 2 over $\mathcal{A}$ for these two modules, it follows easily
that $\mathrm{Ext}^2_{\mathcal{A}}(\mathcal{R},\mathcal{A}(-6))$ equals the full
module of sections of the sheaf
\\$\mathcal{E}xt^2_{\mathcal{O}_{\mathbb{P}^4}}(\tilde{\mathcal{R}},
\mathcal{O}_{\mathbb{P}^4}(-6))$, and $\mathcal{R}$ the full module of sections of
$\tilde{\mathcal{R}}$ (see section 3, lemma 3.3, below, where this argument is made
precise); thus from (1) I infer the isomorphism of $\mathcal{A}-$modules
\begin{equation}
\mathcal{R}=\mathrm{Ext}^2_{\mathcal{A}}(\mathcal{R},\mathcal{A}(-6)),
\end{equation}
which is also an isomorphism of $\mathcal{R}-$modules since it is functorial with
respect to endomorphisms of $\mathcal{R}$ (which follows from the functoriality of the
isomorphisms
$\psi_{\ast}\omega_X=\mathcal{H}om_{\mathcal{O}_Y}(\psi_{\ast}\mathcal{O}_X,\omega_Y)$
and $\mathcal{H}om_{\mathcal{O}_Y}(\psi_{\ast}\mathcal{O}_X,\omega_Y)=\mathcal{E}xt^2_{\mathcal
{O}_{\mathbb{P}^4}}(\psi_{\ast}\mathcal{O}_X,\omega_{\mathbb{P}^4})$ above). The
isomorphism (4) lifts to an isomorphism of minimal graded free resolutions (2) and
(3). In particular,
$\mathrm{rank}\, F_0=\mathrm{rank}\, F_2$, and since
$\mathrm{Ann}_{\mathcal{A}}(\mathcal{R})\neq 0$ one has $\mathrm{rank}\, F_0
-\mathrm{rank}\, F_1+
\mathrm{rank}\, F_2=0$ whence there exists an integer $n$ such that $\mathrm{rank}\,
F_0=\mathrm{rank}\, F_2=n+1$, $\mathrm{rank}\, F_1=2n+2$. For the fact that now (2)
can be symmetrized to give a resolution $\mathbf{R}_{\bullet}$ as in the statement of
the theorem I refer to [Gra], p. 938ff., lemma 2.1 and proposition 2.3., whose proof
applies in the present situation with minor modifications; for a thorough exposition
of the argument that the isomorphism (4) gives rise to a symmetric resolution of
$\mathcal{R}$ cf. also [Wal]. 
\end{proof}
Next I certainly have $p_g(S)\ge 5$ for surfaces $S$ as in definition 1.1, and for
simplicity I assume $p_g(S)=5$ in what follows. As for $K_S^2$ of such surfaces, I
list here:
\begin{itemize}
\item
One can only expect to find canonical surfaces in $\mathbb{P}^4$ with $p_g=5$ and
$q=0$ in the
 range $8\le K^2 \le 54 $. The lower bound follows from Castelnuovo's inequality
$K^2\ge 3p_g+q-7$. The upper bound follows from the Bogomolov-Miyaoka-Yau inequality
$K^2\le 3e(S)$ in combination with Noether's formula $K^2+e(S)=12(1-q+p)$, where
$e(S)$ is the topological Euler characteristic of $S$.
\item
For $K^2=8$ resp. $=9$ the solutions one gets are the complete intersections of type
$(2,4)$ resp. $(3,3)$ (cf. [En], p. 284ff.).
\item
Existence is known in cases $K^2=10,11,12$; the case $K^2=10$ is treated in
[Cil], subsequently also in [Cat 4] (p. 42ff.) and [Ro\ss ] (p. 108ff.), by
approaches different in taste each time. Moreover, in the latter case one has a
satisfactory picture of the moduli space of these surfaces; for $K^2=11,12$ a
partial description of the moduli spaces is in [Ro\ss ].   
\end{itemize}  
Therefore I will also assume $K^2\ge 10$ henceforth.\\
For the case $p_g=5,\: q=0,\: K^2\ge 10$, the numbers $n,\: r_i,\: s_i,\;
i=1,\ldots ,n+1,$ appearing in the resolution $\mathbf{R}_{\bullet}$ of theorem 1.4
are readily calculated; this is done in [Cil], p. 302, prop. 5.3 (cf. also [Cat4],
p. 41, prop. 6.2):
\begin{thm}
For a canonical surface in $\mathbb{P}^4$ with $q=0,\: p_g=5,\: K^2\ge 10$ one has a
resolution of the canonical ring $\mathcal{R}$
\begin{eqnarray}
\begin{CD}\mathbf{R}_{\bullet}:\;0 @>>> \mathcal{A}(-6)\oplus\mathcal{A}(-4)^n@
>{{-\beta^t}\choose{\alpha^t}}>>\mathcal{A}(-3)^{2n+2}\end{CD}\nonumber \\
\begin{CD}@>{(\alpha\,\beta)}>>\mathcal{A}
\oplus\mathcal{A}(-2)^n@>>>\mathcal{R} @>>> 0 ,
\end{CD}
\end{eqnarray}
where $n:=K^2-9$.
\end{thm} 
However, what is important here is that there is a converse to the story told so
far, on which rests the analysis of canonical surfaces done in this work:
\begin{thm}
Let $\mathcal{R}$ be some \underline{algebra} (commutative with $1$) over
$\mathcal{A}=\mathbb{C}[x_0,\ldots ,x_4]$ with minimal graded free resolution as in
(5), with $1\in\mathcal{R}$ corresponding to the first row of $(\alpha\:\beta)$ as
$\mathcal{A}-$module generator. Write
$A:=(\alpha\:\beta)$, $A':=A$ with first row erased, $I_n(A')=$ Fitting ideal of
$n\times n$ minors of $A'$.

Then
$\mathcal{R}$ is a Gorenstein algebra, and if one assumes that
$\mathrm{Ann}_{\mathcal{A}}(\mathcal{R})$ is a prime ideal, then  
$Y:=\mathrm{Supp}(\mathcal{R})\subseteq\mathbb{P}^4$ with its reduced induced
subscheme structure (thus the ideal of polynomials vanishing on $Y$ is
$\mathcal{I}_Y=\mathrm{Ann}_{\mathcal{A}}\,\mathcal{R}$) is an irreducible surface,
and if furthermore one assumes 
$\mathrm{grade}\, I_n(A')\ge 3$ and
$X:=\mathrm{Proj}(\mathcal{R})$ has only rational double points as singularities,
then $X$ is the canonical model of a surface $S$ of general type with $q=0,\: p_g=5,\:
K^2=n+9$. More precisely, writing $\mathcal{A}_Y$ for the homogeneous coordinate ring
of
$Y$, one has that the morphism $\psi : X\to Y\subset\mathbb{P}^4$ induced by the
inclusion
$\mathcal{A}_Y\subset \mathcal{R}$ is a finite birational morphism, and is part of a
diagram\\
\setlength{\unitlength}{1cm}
\begin{center}
\begin{picture}(4,3)
\put(0.7,2.5){$S$}
\put(2.8,2.5){$Y\subset \mathbb{P}^4$}
\put(1.75,0.7){$X$}
\put(1,2.6){\vector(1,0){1.6}}
\put(0.8,2.3){\vector(2,-3){0.8}}
\put(2,1.1){\vector(2,3){0.8}}
\put(1.75,2.7){$\pi$}
\put(2.5,1.5){$\psi$}
\put(1,1.5){$\kappa$} 
\end{picture}\end{center}
where $S$ is the minimal desingularization of $X$, $\kappa$ is the contraction
morphism contracting exactly the (-2)-curves of $S$ to rational double points on $X$,
and the composite $\pi:=\psi\circ\kappa$ is a birational morphism with
$\pi^{\ast}\mathcal{O}_{\mathbb{P}^4}(1)=\mathcal{O}_S(K_S)$ (i.e. is 1-canonical for
$S$).
\end{thm}
\begin{proof}
Taking Hom of (5) into $\mathcal{A}(-6)$ and using the canonical isomorphism from (5)
to its dual, one gets an isomorphism
$\mathcal{R}=\mathrm{Ext}^2_{\mathcal{A}}(\mathcal{R},\mathcal{A}(-6))$, functorial
with respect to endomorphisms of $\mathcal{R}$, which is therefore an isomorphism
of
$\mathcal{R}-$modules. Thus $\mathcal{R}$ is a Gorenstein algebra.\\
Remark that since the ideal of $(n+1)\times (n+1)$ minors of $A$, $I_{n+1}(A)$ (i.e.
the zeroeth Fitting ideal of $\mathcal{R}$), and
$\mathrm{Ann}_{\mathcal{A}}\,\mathcal{R}$ have the same radical, the Eisenbud-Buchsbaum
acyclicity criterion (cf. [Ei], thm. 20.9, p. 500) gives $\mathrm{grade}\,
I_{n+1}(A)=\mathrm{grade}\,\mathrm{Ann}_{\mathcal{A}}\,\mathcal{R}=\mathrm{codim}_{\mathcal{A}}
\,\mathrm{Ann}_{\mathcal{A}}\,\mathcal{R}\ge 2$, whereas also
$\mathrm{grade}\,\mathcal{R}\equiv\mathrm{grade}(\mathrm{Ann}_{\mathcal{A}}\,\mathcal{R}
, 
\mathcal{A})\le\mathrm{projdim}_{\mathcal{A}}\,\mathcal{R}=2$ (cf. e.g. [B-He],
p.25), whence $Y$, defined by the annihilator ideal
$\mathrm{Ann}_{\mathcal{A}}\,\mathcal{R}\subset\mathcal{A}$, is in fact an
irreducible surface.\\  
$\mathcal{R}$ is
CM because
$\mathrm{grade}(\mathcal{R})=\mathrm{projdim}_{\mathcal{A}}(\mathcal{R})$ and thus
$\mathcal{R}$ is a perfect module (cf. [B-He], p. 59, thm. 2.1.5). Next, the
morphism $\psi :X\to Y$ induced by the inclusion $\mathcal{A}_Y\subset \mathcal{R}$
is finite since $\mathcal{R}$ is a finite $\mathcal{A}_Y-$module and thus
$\psi_{\ast}\mathcal{O}_X=\tilde{\mathcal{R}}$ is a finite $\mathcal{O}_Y-$module
over any affine open of $Y$. Now $A'$ is a presentation matrix of
$\mathcal{R}/(\mathcal{A}_Y\cdot 1)$ whence by Fitting's lemma,
$I_n(A')\subset\mathrm{Ann}_{\mathcal{A}}(\mathcal{R}/\mathcal{A}_Y)$ and
$(I_n(A')\cdot \mathcal{A}_Y)\,\mathcal{R}\subset \mathcal{A}_Y$. Since
$(I_n(A')\cdot \mathcal{A}_Y)\neq 0$ ($I_n(A')\not\subset (Y)$ because
$\mathrm{grade}\, I_n(A')\ge 3$) and
$\mathcal{A}_Y$ is an integral domain,
$\exists$ a (homogeneous) nonzerodivisor $d\in \mathcal{A}_Y$ such that $d\cdot
\mathcal{R}\subset
\mathcal{A}_Y$. Now $d$ is also a nonzerodivisor on $\mathcal{R}$ because
$\mathcal{R}$ is a maximal Cohen-Macaulay module over $\mathcal{A}_Y$ (see [Ei], p.
534, prop. 21.9). Thus one gets 
\begin{equation}
\mathcal{R}[d^{-1}]=\mathcal{A}_Y[d^{-1}].
\end{equation}
(Note incidentally that $\mathcal{R}$ is an integral domain because it is contained
in $d^{-1}\mathcal{A}_Y$ and that the algebra structure on $\mathcal{R}$ is uniquely
determined since it is a subalgebra of $\mathcal{A}_Y[d^{-1}]$). From (6) one sees
that $\psi$ gives an isomorphism of function fields $\mathbb{C}(X)=\mathbb{C}(Y)$,
thus is birational.\\
The fact that $X$ has only rational double points as singularities implies that $X$
is locally Gorenstein and the dualizing sheaf $\omega_X$ is invertible,
$\omega_X=\mathcal{O}_X(K_X)$, where $K_X$ is an associated (Cartier) divisor. Now one
can run the argument given in theorem 1.4 in reverse: The sheafified Gorenstein
condition
$\tilde{\mathcal{R}}=\mathcal{E}xt^2_{\mathcal{O}_{\mathbb{P}^4}}(\psi_{\ast}\mathcal{O}_X,
\mathcal{O}_{\mathbb{P}^4}(-6))$ together with $\mathcal{H}om_{\mathcal{O}_Y}(\psi_{\ast}\mathcal{O}_X,\omega_Y)=\mathcal{E}xt^2_{\mathcal
{O}_{\mathbb{P}^4}}(\psi_{\ast}\mathcal{O}_X,\omega_{\mathbb{P}^4})$ and relative
duality for the finite morphism $\psi$ gives
\begin{equation}
\psi_{\ast}\omega_X=\tilde{\mathcal{R}}(1).
\end{equation}
Therefore one finds $\mathcal{R}\cong \bigoplus\limits_{m\ge
0}H^0(X,\mathcal{O}_X(mK_X))$ where the latter is also equal to $\bigoplus\limits_{m\ge
0}H^0(S,\mathcal{O}_S(mK_S))$ with $\kappa : S\to X$ the minimal desingularization of
$X$ and $\kappa$ as in the statement of the theorem. Thus $X$ is the canonical model
of a surface of general type (since $\dim\,\mathcal{R}=3$). Furthermore, since
$\psi^{\ast}\mathcal{O}_{\mathbb{P}^4}(1)=\mathcal{O}_X(K_X)$ and
$\kappa^{\ast}\mathcal{O}_X(K_X)=\mathcal{O}_S(K_S)$, it follows that
$\pi:=\psi\circ\kappa$ is a 1-canonical map for $S$ and clearly a birational morphism,
and
$Y$ is a canonical surface in
$\mathbb{P}^4$.\\
The invariants $p_g(S),\: q(S),\: K_S^2$ are immediately found from the resolution
(5): On the one hand, for the plurigenera one has $P_1=p_g$, $P_m={m \choose
2}K_S^2+\chi (\mathcal{O}_S),\: m\ge 2,$ by Kodaira's formula (cf. [Bom], p. 185),
on the other hand, writing $\mathcal{R}_m$ for the $m$th graded piece of
$\mathcal{R}$, and $\bigoplus\limits_j
\mathcal{A}(-a_{0,j}):=\mathcal{A}\oplus\mathcal{A}(-2)^n,\: \bigoplus\limits_j
\mathcal{A}(-a_{1,j}):=\mathcal{A}(-3)^{2n+2},\: \bigoplus\limits_j
\mathcal{A}(-a_{2,j}):=\mathcal{A}(-6)\oplus\mathcal{A}(-4)^n$, one has
$\dim_{\mathbb{C}}\, \mathcal{R}_m=\sum_{i=0}^2(-1)^i\sum_j {m-a_{i,j}+4 \choose
4}$ from the Hilbert resolution of $\mathcal{R}$ (${k\choose l}=0$ for $k<l$).
Comparing these one finds $p_g=5,\: K^2+6-q=15+n,\: 3K^2+6-q=33+3n$ whence the
invariants are the ones given in the theorem ($q=0$ is also clear since
$\mathcal{R}$ is CM by prop. 1.2).   
\end{proof}
Thus one morally sees that the important question remaining is how the algebra
structure of $\mathcal{R}$ is reflected in the Hilbert resolution resp. to give
necessary and sufficient conditions for the presentation matrix $(\alpha\:\beta)$
such that $\mathcal{R}$ supports the structure of an $\mathcal{A}-$algebra.

\bfseries{Remark 1. }\mdseries
 Let $\mathcal{R}$ be an $\mathcal{A}-$module with resolution (5) such that
$\mathrm{Ann}_{\mathcal{A}}\,\mathcal{R}$ is prime (defining the surface
$Y\subset\mathbb{P}^4$),
$\mathrm{grade}\, I_n(A')\ge 3$, and with a distinguished element, call it 1, in
$\mathcal{R}$ corresponding to the first row of
$(\alpha\:\beta)$. Then it follows from the proof of theorem 1.6 that $\exists$ a
nonzerodivisor $d$ on $\mathcal{R}$ such that $d\cdot
\mathcal{R}\subset\mathcal{A}_Y$; thus if $\mathcal{R}$ is an algebra, it is a
subalgebra of $\mathcal{A}_Y[d^{-1}]$. In particular, in all what follows, an algebra 
structure on $\mathcal{R}$, if it exists, will be \underline{unique}. $\mathcal{R}$
is what is called a finite birational module in [E-U]. In other words, it's a
fractional ideal of $\mathcal{A}_Y$. This entails for example that
$\mathrm{Hom}_{\mathcal{A}_Y}(\mathcal{R},\mathcal{A}_Y)=(\mathcal{A}_Y :_K
\mathcal{R})$ (where $K$ is the quotient field of $\mathcal{A}_Y$) is an ideal of
$\mathcal{A}_Y$ (the so-called conductor of $\mathcal{R}$ into $\mathcal{A}_Y$).

\bfseries{Remark 2. }\mdseries 
With the set-up of theorem 1.6, $V(I_n(A'))=$ nonnormal locus of $Y$. In fact, $\psi
: X\to Y$ is the normalization map, and therefore the sheaf of ideals
$\mathcal{A}nn_{\mathcal{O}_{\mathbb{P}^4}}(\psi_{\ast}\mathcal{O}_X /
\mathcal{O}_Y)=\mathcal{A}nn_{\mathcal{O}_{\mathbb{P}^4}}(\tilde{\mathcal{R}} /
\mathcal{O}_Y)$ defines the nonnormal locus of $Y$. But since $A'$ is a presentation
matrix for $\mathcal{R}/\mathcal{A}_Y$, it is
$\sqrt{\mathrm{Ann}_{\mathcal{A}}(\mathcal{R}/\mathcal{A}_Y)}=\sqrt{I_n(A')}$ (cf.
e.g. [Ei], prop. 20.6, p. 498) and the assertion follows. The assumption
$\mathrm{grade}\, I_n(A')\ge 3$ in the theorem anticipates this fact.

\bfseries{Remark 3. }\mdseries 
Whereas the above theorem is valid without any condition on the singularities of
$Y$, in the sequel it will be sometimes convenient to assume that $Y\subset
\mathbb{P}^4$ has only improper double points as singularities (i.e. points with
tangent cone consisting of two planes spanning $\mathbb{P}^4$); slightly more
generally, the investigation of the ring structure of $\mathcal{R}$ as contained in
theorem 3.1 below can be carried out under the assumption that $Y$ is normal off a
finite number of improper double points. Such $Y$ is sometimes said to have
quasi-ordinary singularities (it is said to have ordinary singularities if it is
smooth off the improper double points). I state here (cf. [Cil], p. 306ff.):
\begin{thm}
Let $\pi : S\to Y\subset \mathbb{P}^4$ be a canonical surface with $q=0,\: p_g=5$. If
$Y$ has ordinary singularities, the number $\delta (Y):={  K_S^2-8 \choose 2}$ is the
number of improper double points of $Y$ (very special case of the "double point
formula of Severi").
\end{thm}

\section[The case $K^2=11$]{Analysis of the case $K^2=11$}
\setcounter{equation}{0}
Let $\pi:S\to Y$ be a canonical surface in
$\mathbb{P}^4$ with
$q=0,\;p_g=5,\;K^2=11$.According to theorem 1.5, one has a resolution
\begin{eqnarray}
\begin{CD}\mathbf{R}_{\bullet}:\;0@>>>\mathcal{A}(-6)\oplus\mathcal{A}(-4)^2@
>{{-\beta^t}\choose{\alpha^t}}>>\mathcal{A}(-3)^6\end{CD}\nonumber\\
\begin{CD}@>{(\alpha\,\beta)}>>\mathcal{A}
\oplus\mathcal{A}(-2)^2@>>>\mathcal{R}@>>>0
\end{CD}
\end{eqnarray}
of the canonical ring $\mathcal{R}$ of $S$. I want to examine how, in
this particular case, $\mathbf{R}_{\bullet}$ encodes the ring structure of
$\mathcal{R}$. In the next section I will treat the case of general
$K^2$, but it may be worthwhile considering $K^2=11$ separately to get
a feeling for the concepts entering the computations and because the
results are accidentally slightly stronger in this case. More notation:
\begin{equation} A:=(\alpha\; \beta)=:\left(
\begin{array}{*{3}{c@{\:}}|*{3}{c@{\:}}} A_1 & A_2 & A_3 & B_1 & B_2 &
B_3 \\ a_1 & a_2 & a_3 & b_1 & b_2 & b_3 \\
a_4 & a_5 & a_6 & b_4 & b_5 & b_6 \end{array}\right),\end{equation}
where the $A_i,\;B_i,\;i=1,2,3,$ are cubic forms, the
$a_j,\;b_j,\;j=1,\ldots,6,$ are linear forms; $A':=A$ with $1^{st}$ row
erased, $I_2(A),I_2(A'):=$Fitting ideals of $2\times 2-$minors
of $A,\;A'$ respectively; $\mathcal{I}_Y:=$ideal of polynomials vanishing on
$Y=\mathrm{Ann}_{\mathcal{A}}\,\mathcal{R}$,
$\mathcal{A}_Y$:=homogeneous coordinate ring of
$Y$. Furthermore I will for simplicity assume that
$Y$ has only improper double points as singularities. Severi's double
point formula then gives that there are 3 of them for $K^2=11$; and we
have: \{3 improper double pt.s of $Y$\}$=V(I_2(A'))$.\\
I claim:
\begin{lemm}Acting on the tableau in (2) with elements
$\left( \begin{array}{*{2}{c@{\:}}}1 & 0\\0 & \varphi\end{array}
\right)$, 
$\varphi \in Gl_2(\mathbb{C})$, from the left, and elements of
$Sp_6(\mathbb{C})$ from the right, one can eventually obtain
\begin{equation}
\tilde{A}=\left( \begin{array}{*{3}{c@{\:}}|*{3}{c@{\:}}}
\tilde{A}_1 & \tilde{A}_2 & \tilde{A}_3 & \tilde{B}_1 & \tilde{B}_2 &
\tilde{B}_3\\
0 & \tilde{a}_2 & \tilde{a}_3 & 0 & \tilde{b}_2 & \tilde{b}_3\\
\tilde{a}_4 & -\tilde{a}_2 & 0 & \tilde{b}_4 & -\tilde{b}_2 & 0
\end{array}\right)=:(\tilde{\alpha}\;\tilde{\beta})
\end{equation}preserving the symmetry:
$\tilde{\alpha}\tilde{\beta}^t=\tilde{\beta}\tilde{\alpha}^t$. The
$\tilde{a_i},\;
\tilde{b_i}$ are linear forms s.t.
\begin{eqnarray*}
V(\tilde{a}_2,\tilde{a}_3,\tilde{b}_2,\tilde{b}_3)=\{1^{st}\; improper\;
dbl.\; pt.\}\\
V(\tilde{a}_4,\tilde{a}_2,\tilde{b}_4,\tilde{b}_2)=\{2^{nd}\; improper\;
dbl.\; pt.\}\\
V(\tilde{a}_4,\tilde{a}_3,\tilde{b}_4,\tilde{b}_3)=\{3^{rd}\; improper\;
dbl.\; pt.\}.
\end{eqnarray*}
(The $\tilde{A}_j,\;\tilde{B}_j,\;j\in\{1,2,3\}$ are of course
cubics, linear combinations of the $A_j,\;B_j$).    
\end{lemm}
\begin{proof}Write $\alpha_{\mu}$ resp. $\beta_{\nu}$ for the $\mu -$th
resp. $\nu -$th column of $\alpha$ resp. $\beta$. Let's make a list of
some useful allowable operations on $A$, i.e. such that they preserve the
symmetry:
\begin{itemize}
\item[(i)] Elementary operations on rows: indeed,
$\forall g=\left( \begin{array}{*{2}{c@{\:}}}1 & 0\\0 & \varphi\end{array}
\right),\: 
\varphi \in Gl_2(\mathbb{C}):\:\alpha\beta^t=\beta\alpha^t\Rightarrow (g\alpha
)(g\beta)^t=(g\beta)(g\alpha)^t$
\item[(ii)]
For $\lambda\in\mathbb{C}$ and $\mu$ a fixed but arbitrary column index, adding
$\lambda\beta_{\mu}$ to $\alpha_{\mu}$: Remark that both sides of each of the
equations $\sum_i \alpha_{hi}\beta_{li}=\sum_i \beta_{hi}\alpha_{li}$ are just changed
by a summand $\lambda\beta_{h\mu}\beta_{l\mu}$. This operation is of course as well
applicable with the r\^{o}les of $\alpha$ and $\beta$ interchanged.
\item[(iii)]
For $\lambda\in\mathbb{C}$ and $\mu ,\:\nu$ column indices, adding
$\lambda\beta_{\nu}$ to $\alpha_{\mu}$ and at the same time adding
$\lambda\beta_{\mu}$ to $\alpha_{\nu}$: Both sides of each of the equations $\sum_i
\alpha_{hi}\beta_{li}=\sum_i \beta_{hi}\alpha_{li}$ change by a summand $\lambda
(\beta_{h\nu}\beta_{l\mu}+\beta_{h\mu}\beta{l\nu})$. [(ii) is thus a special case of
(iii) with $\mu =\nu$]; the same operation also with the r\^{o}les of $\alpha,\:\beta$
interchanged.
\item[(iv)]
For $\lambda\in\mathbb{C}$ and $\mu\neq\nu$ column indices, adding
$\lambda\alpha_{\nu}$ to $\alpha_{\mu}$ and simultaneously subtracting
$\lambda\beta_{\mu}$ from $\beta_{\nu}$; this is O.K. since it corresponds to changing
the left side of $\sum_i \alpha_{hi}\beta_{li}=\sum_i \beta_{hi}\alpha_{li}$ by a
summand
$\lambda (\alpha_{h\nu}\beta_{l\mu}-\alpha_{h\nu}\beta_{l\mu})=0$, and the right side
by a summand $\lambda(\beta_{h\mu}\alpha_{l\nu}-\beta_{h\mu}\alpha_{l\nu})=0$; the same
operation also with the r\^{o}les of
$\alpha,\:\beta$ interchanged.
\item[(v)] For $\mu\neq\nu$, interchanging columns $\alpha_{\mu},\:\alpha_{\nu}$ and
at the same time interchanging columns $\beta_{\mu},\:\beta_{\nu}$, which clearly
preserves the symmetry.
\item[(vi)] For a column index $\mu$, multiplying column $\alpha_{\mu}$ by $(-1)$ and
then interchanging columns $-\alpha_{\mu}$ and $\beta_{\mu}$ (i.e. the substitution
$\alpha_{\mu}\mapsto\beta_{\mu},\:\beta_{\mu}\mapsto -\alpha_{\mu}$): Namely, $\sum_i
\alpha_{hi}\beta_{li}=\sum_i \beta_{hi}\alpha_{li}\Leftrightarrow \sum_{i\neq\mu}
\alpha_{hi}\beta_{li}-\beta_{h\mu}\alpha_{l\mu}=\sum_{i\neq\mu}
\beta_{hi}\alpha_{li}-\alpha_{h\mu}\beta_{l\mu}$.
\end{itemize}
Call these operations (Op). Remark that (Op), (ii)-(vi)
correspond to multiplication on $A$ from the right by symplectic $6\times 6$
matrices. In fact, more systematically, one sees that since symplectic matrices
$\left(
\begin{array}{cc} S_1 & S_2\\ S_3 & S_4 \end{array}\right) \in Gl_{2n+2}(\mathbb{C})$,
$S_1,\: S_2,\: S_3,\: S_4$ $(n+1)\times (n+1)$ matrices, can be characterized by
equations $S_1S_2^t=S_2S_1^t,\: S_3S_4^t=S_4S_3^t,\: S_1S_4^t-S_2S_3^t=I_{n+1}$, one
has for $A=(\alpha\:\beta )$ an $(n+1)\times (n+1)$ matrix with $\alpha\beta^t$
symmetric (as in thm. 1.5) that also $(\alpha S_1+\beta S_3)(\alpha S_2+\beta S_4)^t$
is symmetric (this is also immediate because the symmetry condition can be rephrased
as saying that, for each choice of homogeneous coordinate vector $(x_0:\ldots : x_4)$
in $\mathbb{P}^4$, the rows of $(\alpha\:\beta)$ span an isotropic subspace for the
standard symplectic form on $\mathbb{C}^{2n+2}$, and a matrix is symplectic iff its
transpose is).\\

First a general remark: Given a matrix of linear forms, call a generalized
row of this matrix an arbitrary linear combination of the rows with not all
coefficients zero. Then the locus where the rows are linearly dependent is the union,
over all generalized rows, of the linear spaces cut out by the linear forms which are
the entries of the generalized row.  Therefore I can assume
$A'=\left(\begin{array}{*{3}{c@{\:}}|*{3}{c@{\:}}}
a_1 & a_2 & a_3 & b_1 & b_2 & b_3\\
a_4 & a_5 & a_6 & b_4 & b_5 & b_6
\end{array}
\right)$ to be such that one of the improper double points is given by the
vanishing of the linear forms in the upper row of $A'$, the second one by
the vanishing of the linear forms in the lower row of $A'$, and the third
as the zero set of the linear forms gotten by adding the two rows
together.\\ The rest of the proof is a game on the tableau $A'$, using
(Op) and the symmetry, and deriving Koszul sequences from the fact that
the rows of
$A'$ resp. their sum define 3 distinct points. To ease notation, I will
treat the $a_i,\;b_i,\;i=1,\ldots,6$, and $A'$ as dynamical
variables. For clarity's sake, I will box certain assumptions in the course
of the following argument, especially when they are cumulative.
\\Using (Op),(v)/(vi), then (iv) and finally (iii) one gets
\begin{equation}A'=\left(\begin{array}{*{3}{c@{\:}}|*{3}{c@{\:}}} 0 & a_2 & a_3 & b_1
& b_2 & b_3\\ a_4 & a_5 & a_6 & b_4 & b_5 & b_6
\end{array}
\right);\end{equation} \fbox{$a_4=0$}: Then
$D:=\left\{\mathrm{rk}\left( \begin{array}{*{5}{c@{\:}}} a_2 & a_3 & b_1 &
b_2 & b_3\\a_5 & a_6 & b_4 & b_5 &b_6 
\end{array}
\right)\le 1
\right\}=\{\mathrm{3\; dbl.\; pt.s}\}.$ But the determinantal locus $D$
has the expected codimension since the generic $2\times
5-$matrix degenerates in codimension 4. Therefore, by
Porteous' formula (cf. [A-C-G-H], p. 90ff.), its degree is also the expected one,
namely 5, a contradiction. Therefore quite generally the
possibility of a zero column is excluded.\\\fbox{$a_4\neq 0$}: Use (Op),
(v), (vi), (iv) , (iii) in this order to put a zero in place of
$b_6$ ($a_4,\; a_5,\;a_6,\;b_5,\;b_6$ are dependent!) :
\begin{equation}A'=\left(\begin{array}{*{3}{c@{\:}}|*{3}{c@{\:}}} 0 & a_2
& a_3 & b_1 & b_2 & b_3\\ a_4 & a_5 & a_6 & b_4 & b_5 & 0
\end{array}
\right).\end{equation} I claim that now \fbox{$a_2,\;b_1,\;b_2,\;b_3$ are
dependent}: For if they are independent I can also assume
$a_4,\;a_5,\;a_6,\;b_5$ independent (otherwise interchange rows and use
(Op), (v) and (vi)).Symmetry gives:
$b_1a_4+b_2a_5+b_3a_6+b_5\cdot(-a_2)=0$, which is a Koszul relation
saying $\exists$ antisymmetric matrices of scalars $S,\; \tilde{S}$
such that
\[\left(\begin{array}{c}a_4\\a_5\\a_6\\b_5\end{array}\right)=\tilde{S}\left(
\begin{array}{c}b_1\\b_2\\b_3\\-a_2\end{array}\right),\;\left(\begin{array}{c}b_1\\b_2\\b_3\\-a_2\end{array}\right)=S\left(
\begin{array}{c}a_4\\a_5\\a_6\\b_5\end{array}\right),\]$\tilde{S}S=I$, $S,\;\tilde{S}$
are invertible. Now interchange the $4^{th}$ and $5^{th}$ columns of $A'$
and multiply by $\left( \begin{array}{c|cc}S^t & \cdots  & 0\\ \hline
\vdots & 1&0\\ 0&0&1\end{array}
\right)$ on the right. (This will in general destroy the symmetry but
preserve the points that are defined by the rows of $A'$ and their
sum; this operation is only used to derive a contradiction). The
second row of the transformed matrix  is then
$(b_1,b_2,b_3,-a_2,b_4,0)$, and one sees that it either defines
$\emptyset$ or the same point as the first row, a contradiction
because I assumed the points defined by the rows of $A'$ to be
distinct. Therefore $a_2,\;b_1,\;b_2,\;b_3$ are dependent.\\I claim
further that then \fbox{$a_2,\;b_2,\;b_3$  are independent}. Suppose
not.  Since the possibility of a zero column was excluded for reasons of
degree above, I can then use (Op), (iii) and if necessary (vi) to get $A'=\left(\begin{array}{*{3}{c@{\:}}|*{3}{c@{\:}}} 0 &
0 & a_3 & b_1 & b_2 & b_3\\ a_4 & a_5 & a_6 & b_4 & b_5 & 0
\end{array}
\right)$. Here  $a_3,\;b_1,\;b_2,\;b_3$ are independent. I have 2 cases:
\begin{enumerate}
\item
$a_4,\;a_5,\;a_6$ are independent. Then the symmetry gives that $\exists$
antisymmetric matrices of scalars $T,\;\tilde{T}$ s.t.
$\left(\begin{array}{c}a_4\\a_5\\a_6\end{array}\right)=T\left(
\begin{array}{c}b_1\\b_2\\b_3\end{array}\right),\;\left(\begin{array}{c}b_1\\b_2\\b_3\end{array}\right)=\tilde{T}\left(
\begin{array}{c}a_4\\a_5\\a_6\end{array}\right);$ but then
$T,\;\tilde{T}$ are invertible contradicting the fact that a
skewsymmetric matrix of odd-dimensional format has determinant zero.
\item
$a_4,\;a_5,\;a_6$ are dependent. Since no zero column can occur, I can
use (Op), (iv) to write $A'=\left(\begin{array}{*{3}{c@{\:}}|*{3}{c@{\:}}} 0 &
0 & a_3 & b_1 & b_2 & b_3\\ a_4 & a_5 & 0 & b_4 & b_5 & 0
\end{array}
\right)$; but the symmetry $a_4b_1=-a_5b_2$ tells me I am left with
discussing the case $A'=\left(\begin{array}{*{3}{c@{\:}}|*{3}{c@{\:}}} 0 &
0 & a_3 & -a_5 & a_4 & b_3\\ a_4 & a_5 & 0 & b_4 & b_5 & 0
\end{array}
\right)$. But then the points defined by the second row and the sum of
the rows coincide, or the linear forms in the sum of the rows define
$\emptyset$, a contradiction. 
\end{enumerate}Using the last two boxed assumptions and (Op), (iii) and
then (iv), I can pass from the shape of $A'$ in (5) to \begin{equation}A'=\left(\begin{array}{*{3}{c@{\:}}|*{3}{c@{\:}}} 0 & a_2
& a_3 & 0 & b_2 & b_3\\ a_4 & a_5 & a_6 & b_4 & b_5 & 0
\end{array}
\right).\end{equation} Now I play the game again, but this time it is
quicker. I claim:\\ \fbox{$a_5,\;a_6,\;b_5$ dependent}: If not, the
symmetry $a_5b_2+a_6b_3+b_5(-a_2)=0$ gives as above the existence of
$3\times 3$ invertible skew-symmetric matrices, a contradiction. But I
also claim: \fbox{$a_5,\;b_5$ are independent}: Otherwise I get, using (Op), (ii) and
possibly (vi), $A'=\left(\begin{array}{*{3}{c@{\:}}|*{3}{c@{\:}}} 0 & a_2
& a_3 & 0 & b_2 & b_3\\ a_4 & 0 & a_6 & b_4 & b_5 & 0
\end{array}
\right),$ and using  the symmetry $a_2b_5=a_6b_3$, I must look at  $A'=\left(\begin{array}{*{3}{c@{\:}}|*{3}{c@{\:}}} 0 &
a_6 & a_3 & 0 & b_2 & b_5\\ a_4 & 0 & a_6 & b_4 & b_5 & 0
\end{array}
\right).$ But then either the points defined by the second row of
$A'$ and the sum of its rows resp. coincide, or the linear forms in
the sum of the rows define $\emptyset$, a contradiction. 
Using the previous 2 boxed assumptions and (Op), (iv) and then
(iii), I can pass from (6) to \begin{equation}A'=\left(\begin{array}{*{3}{c@{\:}}|*{3}{c@{\:}}} 0 & a_2
& a_3 & 0 & b_2 & b_3\\ a_4 & a_5 & 0 & b_4 & b_5 & 0
\end{array}
\right).\end{equation} Invoking the symmetry $a_2b_5=a_5b_2$ a last
time, I am through:
\[A'=\left(\begin{array}{*{3}{c@{\:}}|*{3}{c@{\:}}} 0 & a_2
& a_3 & 0 & b_2 & b_3\\ a_4 & - a_2 & 0 & b_4 & - b_2 & 0
\end{array}
\right).\]
   
\end{proof}
The above lemma implies that for $K^2=11$ the one half of the
"Eisenbud-Ulrich ring condition" on the presentation matrix
$A=(\alpha\;\beta)$ of $\mathcal{R}$ is more or less automatical:
\begin{prop}Assume $A=(\alpha\;\beta)$ with $\alpha\beta^t=\beta\alpha^t$
is as in (1) (precisely, I assume A is a 3$\times6-$matrix with
first row cubic forms and other rows linear forms, and the locus where
$A$ drops rank consists of 3 distinct points (the improper double points
of Y)). Then \begin{equation} I_2(A)=I_2(A').\end{equation}
\end{prop}
{\bfseries Remark.} The theorem of Eisenbud \&
Ulrich says, in this special case: If (8) holds  and in addition
grade$(I_2(A'))\ge5$ or
$I_2(A')$ is radical, then this suffices to give
and determine a ring structure on $\mathcal{R}^{\ast\ast}$ (double dual
with respect to the
$\mathcal{A}_Y-$module structure of
$\mathcal{R}$). I will show below (cf. lemma 3.3) that in the situation I am in
$\tilde{\mathcal{R}}=\tilde{\mathcal{R}}^{\ast\ast}$ as sheaves, but one remarks that
$I_2(A')$ is generated by quadratics and the 3 points it defines trivially lie in a
hyperplane in $\mathbb{P}^4$, so $I_2(A')$ is not radical; nor is
grade$(I_2(A'))\ge5$ since here $I_2(A')$ defines a codimension 4
subset. Cf. section 3 below that their theorem does not apply in case of
higher $K^2$ either. 
\\
\begin{proof}The idea, if there is any, is simply that I can write $A$ as
in Lemma 2.1, and then for each of the 3 points in the degeneracy locus
of $A$ I get a bunch of Koszul relations which allow me to check (8) by
explicit computation. So write:
\[A=\left(\begin{array}{*{3}{c@{\:}}|*{3}{c@{\:}}}A_1&A_2&A_3&B_1&B_2&B_3
\\0 & a_2 & a_3 & 0 & b_2 & b_3\\ a_4 & - a_2 & 0 & b_4 & - b_2 & 0
\end{array}
\right),\]$v_1:=(a_2,\;a_3,\;b_2,\;b_3)^t,\;v_2:=(a_4,\;-
a_2,\;b_4,\;- b_2)^t,\;v_3=(a_4,\;a_3,\; b_4,\;b_3)^t$ regular sequences. Putting
$W_1:=(-B_2,\;-B_3,\:A_2,\;A_3)^t,\:W_2:=(-B_1,\;-B_2,$ \\$A_1,\;A_2)^t,\;W_3:=(B_1,\;
B_3,\;-A_1,\;-A_3)^t$,
the symmetry amounts to:\[W_1^tv_1=0,\;W_2^tv_2=0,\;W_3^tv_3=0,\]where
the last equation is obtained by adding the first two together.
These are Koszul relations saying that  $\exists$ skew-symmetric
matrices $P,\;Q,\;R$ of quadratic forms such that
\begin{equation}W_1=Pv_1,\;W_2=Qv_2,\;W_3=Rv_3.
\end{equation}$I_2(A')$ is generated by all possible products of
elements of the first and second row of$A'$ resp. except
$a_2^2,\;a_2b_2,\;b_2^2$. A direct computation using the relations
(9) shows $I_2(A)\subseteq I_2(A')$. (To take up an example, look  at
$A_2b_2-B_2a_2\in I_2(A)$. Writing out the second resp. fourth
vector component of the equation $W_2=Qv_2$ in (9) gives
\begin{eqnarray*}-B_2=Q_{21}a_4+Q_{23}b_4-Q_{24}b_2 \\
A_2=Q_{41}a_4+Q_{24}a_2+Q_{43}b_4\end{eqnarray*}Multiplying the first by
$a_2$, the second by $b_2$ and adding gives $A_2b_2-B_2a_2\in I_2(A')$
since $a_2a_4,\;a_2b_4,\;b_2a_4,\;b_2b_4\in I_2(A')$. Here the
skew-symmetry of $Q$ is relevant. Similarly for the other minors.)    
\end{proof}
In section 3 (cf. p. 25ff. below) I will show that the condition $I_2(A')=I_2(A)$ of
proposition 2.2 (or slightly weaker $\overline{I_2(A)}=\overline{I_2(A')}$ where the
bar denotes saturation) implies under the assumption that $Y$ has only improper double
points as singularities that a ring structure is given to and determined on 
$\mathcal{R}$ via
$\mathcal{R}=\mathrm{Hom}_{\mathcal{A}_Y}(\overline{I_2(A')}\cdot\mathcal{A}_Y,\overline{I_2(A')}\cdot\mathcal{A}_Y).$
 I'm sorry I have to refer forward to this result but I don't know any simpler proof
in the special case $K^2=11$ than the one that applies uniformly to all higher $K^2$
as well. Assuming this for the moment and combining what has been said so far with
theorem 1.6, what one gets out of the above discussion is this: \\
The datum (D) of
\itshape
\begin{quote}
a matrix $A=\left(\begin{array}{*{3}{c@{\:}}|*{3}{c@{\:}}}A_1&A_2&A_3&B_1&B_2&B_3
\\0 & a_2 & a_3 & 0 & b_2 & b_3\\ a_4 & -a_2 & 0 & b_4 & -b_2 & 0
\end{array}
\right)$ with the $A_i,\: B_i,\: i=1,\ldots , 3$ cubic forms, $a_2,\: a_3,\: a_4,\:
b_2,\: b_3,\: b_4$ linear forms on $\mathbb{P}^4$ satisfying the symmetry $A_2 b_2+A_3
b_3+B_2(-a_2)+B_3(-a_3)=0$, $A_1b_4+A_2(-b_2)+B_1(-a_4)+B_2a_2=0$, plus the open
conditions that, with $\mathcal{R}:=\mathrm{coker}\, A$,
$\mathrm{Ann}_{\mathcal{A}}\,\mathcal{R}$ be prime and 
$Y=\mathrm{Supp}(\mathcal{R})\subset \mathbb{P}^4$ (with its reduced induced closed 
subscheme structure) be an (irreducible) surface with only singularities $3$ improper
double points given by
$V(a_2,a_3,b_2,b_3)$, $V(a_4,-a_2,b_4,-b_2)$ and $V(a_4,a_2,b_4,b_2)$ and
$X=\mathrm{Proj}\,\mathcal{R}$ have only rational double points as singularities,
modulo graded automorphisms of
$\mathcal{A}\oplus\mathcal{A}(-2)^2$ resp.
$\mathcal{A}(-3)^6$ (acting on $A$ from the right resp. left) which preserve the
normal form of $A$ just described, modulo automorphisms of $\mathbb{P}^4$,
\end{quote}
\upshape
is equivalent to the datum (D') of
\itshape
\begin{quote}
a canonical surface $\pi :S\to Y\subset \mathbb{P}^4$ with $q=0,\: p_g=5,\: K^2=11$
such that $Y$ has only improper double points as singularities, modulo isomorphism.
\end{quote}
\upshape 
The benefit of the normal form in (D) to which the presentation matrices
$A=(\alpha\:\beta)$ of the canonical rings $\mathcal{R}$ of the afore-mentioned
surfaces can be reduced is that the symmetry condition ($\approx$ ring condition)
$\alpha\beta^t=\beta\alpha^t$ can be explicitly solved (for $K^2>11$ I find no such
normal form). Furthermore this can be used to describe the set of isomorphism classes
of surfaces in (D') inside their moduli space $\mathfrak{M}_{K^2,\chi
}=\mathfrak{M}_{11,6}$:\\ 
First one notes that as in the proof of prop. 2.2, the symmetry condition in (D)
amounts to the existence of skew-symmetric $4\times 4$ matrices $P=(P_{ij})$ and
$Q=(Q_{ij})$ of quadratic forms such that
$(-B_2,-B_3,A_2,A_3)^t=P(a_2,a_3,b_2,b_3)^t$ and
$(-B_1,-B_2,A_1,A_2)^t=Q(a_4,-a_2,b_4,-b_2)^t$. Of course there is some ambiguity in
choice of the $(P_{ij}),\: (Q_{ij})$, for the Koszul complexes
$\mathbf{K}_{\bullet}(a_2,a_3,b_2,b_3)$ and $\mathbf{K}_{\bullet}(a_4,-a_2,b_4,-b_2)$
associated to these regular sequences
\begin{eqnarray*}
 {A}(-4)\stackrel{d_3}{\rightarrowtail}\mathcal{A}(-3)^4
\stackrel{d_2}{\longrightarrow}\mathcal{A}(-2)^6\stackrel{d_1}{\longrightarrow}\mathcal{A}(-1)^4
\stackrel{d_0}{\longrightarrow}\mathcal{A}\twoheadrightarrow \mathcal{A}
/(a_2,a_3,b_2,b_3),\\ 
 {A}(-4)\stackrel{d_3'}{\rightarrowtail}\mathcal{A}(-3)^4
\stackrel{d_2'}{\longrightarrow}\mathcal{A}(-2)^6\stackrel{d_1'}{\longrightarrow}\mathcal{A}(-1)^4
\stackrel{d_0'}{\longrightarrow}\mathcal{A}\twoheadrightarrow \mathcal{A}
/(a_4,-a_2,b_4,-b_2)
\end{eqnarray*}
show that e.g. the vector $(P_{ij})_{i<j}\in \mathcal{A}(-2)^6_4$ is only determined
up to addition of $d_2(\underline{l})$ where $\underline{l}\in \mathcal{A}(-3)^4_4$
is a vector of linear forms, and two $\underline{l}$'s give rise to the same
$(P_{ij})_{i<j}$ iff they differ by $d_3(s)$ where $s\in \mathcal{A}(-4)_4$ is a
complex scalar. In other words, $\dim_{\mathbb{C}}(\ker(d_1)_4)=19$ and effectively,
instead of the $(P_{ij})_{i<j}$, one chooses $(\overline{P}_{ij})_{i<j}\in
\mathcal{A}(-2)^6_4 / d_2(\mathcal{A}(-3)^4_4/d_3(\mathcal{A}(-4)_4))$. Similarly for
the $(Q_{ij})$.\\
Next it is clear that whereas now \fbox{$P_{24}$} and \fbox{$Q_{13}$} are subject to no
further relations, for the $\{ (P_{ij})_{i<j} \} -\{ P_{24}\}$ and $\{ (Q_{ij})_{i<j}
\} -\{ Q_{13}
\}$ the relations
\begin{eqnarray*}
A_2=-P_{13}a_2-P_{23}a_3+P_{34}b_3 \hspace{1cm}   B_2=-P_{12}a_3-P_{13}b_2-P_{14}b_3   
\\ A_2=-Q_{14}a_4+Q_{24}a_2-Q_{34}b_4 \hspace{1cm}  B_2=Q_{12}a_4-Q_{23}b_4+Q_{24}b_2
\end{eqnarray*}
imply relations
\begin{eqnarray}
Q_{14}a_4+(P_{13}+Q_{24})(-a_2)+(-P_{23})a_3+Q_{34}b_4+P_{34}b_3=0\\
Q_{12}a_4+P_{12}a_3+(-Q_{23})b_4+(P_{13}+Q_{24})b_2+P_{14}b_3=0.
\end{eqnarray}
I claim that I can assume that the sequences $(a_4,-a_2,a_3,b_4,b_3)$ and\\
$(a_4,a_3,b_4,b_2,b_3)$ are both regular whence (10) and (11) would be Koszul
relations. According to the normal form of the matrix $A$ given in (D), 
$a_4,a_3,b_4,b_3$ are independent (and define one of the improper double points of
$Y$). Assume both $-a_2$ and $b_2$ were expressible in terms of the latter. Then
$V(a_2,a_3,b_2,b_3)$ and $V(a_4,a_3,b_4,b_3)$ would not give distinct points,
contradiction. Therefore at least one of the sequences $(a_4,-a_2,a_3,b_4,b_3)$ and
$(a_4,a_3,b_4,b_2,b_3)$ is regular. But if one of them, $(a_4,-a_2,a_3,b_4,b_3)$ say,
is not regular, then replacing $a_2$ with $a_2+b_2$ (which corresponds to applying
once (Op), (ii) to the matrix $A$) the sequence $(a_4,-(a_2+b_2),a_3,b_4,b_3)$ will
be regular. Similarly if $(a_4,a_3,b_4,b_2,b_3)$ fails to be regular. \\
Therefore considering (10) and (11) as Koszul relations, one gets two skew-symmetric
$5\times 5$ matrices $L=(L_{kl})$ and $M=(M_{kl})$ of linear forms such that
\begin{eqnarray}   
(Q_{14},P_{13}+Q_{24},-P_{23},Q_{34},P_{34})^t=L(a_4,-a_2,a_3,b_4,b_3)^t\\
(Q_{12},P_{12},-Q_{23},P_{13}+Q_{24},P_{14})^t=M(a_4,a_3,b_4,b_2,b_3)^t.
\end{eqnarray}
Call $\underline{a}:=(a_4,-a_2,a_3,b_4,b_3)$, $\underline{a}':=(a_4,a_3,b_4,b_2,b_3)$.
Again looking at Koszul complexes
\begin{eqnarray*}
\mathcal{A}(-5)\stackrel{D_4}{\rightarrowtail}
{A}(-4)^5\stackrel{D_3}{\longrightarrow}\mathcal{A}(-3)^{10}
\stackrel{D_2}{\longrightarrow}\mathcal{A}(-2)^{10}\stackrel{D_1}{\longrightarrow}\mathcal{A}(-1)^5
\stackrel{D_0}{\longrightarrow}\mathcal{A}\twoheadrightarrow \mathcal{A}
/\underline{a},\\ 
\mathcal{A}(-5)\stackrel{D_4'}{\rightarrowtail}
{A}(-4)^5\stackrel{D_3'}{\longrightarrow}\mathcal{A}(-3)^{10}
\stackrel{D_2'}{\longrightarrow}\mathcal{A}(-2)^{10}\stackrel{D_1'}{\longrightarrow}\mathcal{A}(-1)^5
\stackrel{D_0'}{\longrightarrow}\mathcal{A}\twoheadrightarrow \mathcal{A}
/\underline{a}' 
\end{eqnarray*}
one sees that whereas e.g. the $(L_{kl})$ are not unique, the
$(\overline{L}_{kl})_{k<l}\in$\\ $\mathcal{A}(-2)^{10}_3/ D_2(\mathcal{A}(-3)^{10}_3)$
are, and $\dim_{\mathbb{C}}\, (\ker (D_1))_3=10$. Likewise for the $(M_{kl})$.\\
Now equations (12) and (13) should be interpreted as saying that after one of $P_{13}$
and $Q_{24}$, \fbox{$P_{13}$} say, is chosen freely, the other $P$'s and $Q$'s in
(12) and (13) are determined by $L,\: M,\: \underline{a},\: \underline{a}'$.\\
Furthermore one remarks that then the 6 \fbox{$(L_{kl})_{k<l,k\neq 2,l\neq 2}$} and
the 6\\ \fbox{$(M_{kl})_{k<l,k\neq 4,l\neq 4}$} satisfy no further relations, but the
other ones enter in the following relation resulting from equating the 2nd resp. 4th
vector components of (12) resp. (13):
\begin{equation}
(M_{14}-L_{12})a_4+(M_{24}+L_{23})a_3+(M_{34}+L_{24})b_4+(-M_{45}+L_{25})b_3=0.
\end{equation}
The sequence $(a_4,a_3,b_4,b_3)$ is regular by the characterization of the normal form
of $A$ given in (D). One therefore infers the existence of a $4\times 4$
skew-symmetric matrix $S=(S_{rs})$ of complex scalars such that
\begin{equation}
(M_{14}-L_{12},M_{24}+L_{23},M_{34}+L_{24},-M_{45}+L_{25})^t=S(a_4,a_3,b_4,b_3)^t
\end{equation} 
and one notes that the $(S_{rs})$ are then uniquely determined from equation (14).
Moreover upon choosing \fbox{$M_{14},\: M_{24},\: M_{34}, \: M_{45}$} arbitrarily, I
can recover $L_{12},\: L_{23},\: L_{24},\: L_{25}$ from $S$ and $(a_4,a_3,b_4,b_3)$
using (15); and the 6 scalars \fbox{$(S_{rs})_{r<s}$} are not subject to any other
relation in the present set-up.

To get back to the study of the moduli space of surfaces in (D'), fit together the
$a_t,\: b_t,\: t=2,\ldots,4$ and all the boxed objects above into one big affine space
of parameters: 
\begin{displaymath}
\mathcal{P}=\left\{ \begin{array}{l|l}
  & t \in \{ 2,3,4\};\: 
k,\: l,\:  \kappa ,\: \lambda \in \{1,\ldots ,5\},\: \kappa < \lambda ,
 \\ (a_t,b_t,P_{24},Q_{13}, &  k<l,\: k\neq 2,\: l\neq 2;\:  r,\: s\in
\{1,\ldots ,4\},\: r<s;
\\
 P_{13},L_{kl},M_{\kappa \lambda},S_{rs}) & \mathrm{and}\: P_{24},\:Q_{13},\:P_{13}\:
\mathrm{quadratic}, \: a_t,\: b_t,\:  L_{kl},  \\ & M_{\kappa
\lambda}\: \mathrm{linear}\:  \mathrm{in}\: \mathrm{the}\:
\mathrm{hom.}\:\mathrm{coord.}\: (x_0:\ldots :x_4), \\  &
S_{rs} \: \mathrm{complex}\: \mathrm{scalars}. 
\end{array}
\right\}.
\end{displaymath}
Counting one finds that there are 3 quadratic forms, 22 linear forms and 6 scalars in
$\mathcal{P}$, depending on 45, 110 and 6 parameters respectively, whence I have
$\mathcal{P}=\mathbb{A}^{161}$. \\
According to the above discussion, for each choice in an open set of $\mathcal{P}$ one
gets a matrix $A$ meeting the requirements in (D) and a ring $\mathcal{R}$ which is
the canonical ring of a surface of general type $S$ as in (D'). In other words the
parameter space for the canonical rings of the surfaces in (D') is a projection of an
open set of $\mathcal{P}$. (In fact it would be necessary to show that this open set
is non-empty; this is possible, making general choices in $\mathcal{P}$ and verifying
that one gets a matrix $A$ fulfilling the open conditions in (D) e.g. with the help of
a computer algebra package like MACAULAY; anyway, the existence of surfaces in (D')
has been established by Ro\ss berg, cf. [Ro\ss ], 112ff., whence I do not carry out
what I said before).\\
In particular, by the preceding remark one finds that the surfaces in (D') form an
irreducible open set $\mathfrak{U}$ inside their moduli space, and $\mathfrak{U}$ is
unirational (since $\mathcal{P}$ is rational). \\
To calculate the dimension of $\mathfrak{U}$ I note that I have 3 groups acting on the
set of normal forms of matrices $A$ in (D):
\begin{enumerate}
\item
$G=PGl(5)=Aut(\mathbb{P}^4)$ changing coordinates $(x_0:\ldots :x_4)\mapsto
(y_0:\ldots : y_4)$, $\dim \, G=24$.
\item
$ H=\left\{ \mathrm{graded}\: \mathrm{auto.'s}\: 
\mathrm{of}\:
\mathcal{A}\oplus\mathcal{A}(-2)^2 \: \mathrm{of}\:\mathrm{the}\:
\mathrm{form}
\left(
\begin{array}{ccc} s_1 & q_1 & q_2 \\
0   & s_2 & 0   \\
0   &  0  & s_2 
\end{array}
\right)   
\right\}$\\ with $s_1,\: s_2\in\mathbb{C}\backslash \{ 0\}$ and $q_1,\: q_2$
quadratic. Here
$\dim
\, H=32$.
\item
$L=\left\{ \mathrm{group}\:\mathrm{of}\:\mathrm{matrices}\: \small
\left( \begin{array}{cccccc}
\lambda_1 & 0  & 0 & \mu_1 & 0 & 0\\
0 & \lambda_2 & 0 & 0 & \mu_2 & 0 \\
0 & 0 & \lambda_3 & 0 & 0 & \mu_3 \\
\mu_4 & 0 & 0 & \lambda_4 & 0 & 0\\
0 & \mu_5 & 0 & 0 & \lambda_5 & 0 \\
0 & 0 & \mu_6 & 0 & 0 & \lambda_6\\
\end{array}\right)  \right\}$ \normalsize
$\cap Sp_6(\mathbb{C})$, \\
where $\lambda_i\in\mathbb{C}\backslash \{ 0 \},\: \mu_i\in\mathbb{C}$, $i=1,\ldots
,6$, and $Sp_6(\mathbb{C})$ denotes the group of symplectic $6\times 6$ matrices. I
have $\dim\, L=9$.
 
\end{enumerate}
Thus one can calculate an upper bound for the dimension of $\mathfrak{U}$ as follows:
\begin{eqnarray*}
161(\dim  \, \mathcal{P})-38(\dim _{\mathbb{C}}(\ker (d_1)_4)+\dim _{\mathbb{C}}(\ker 
(d_1')_4))-20
(\dim _{\mathbb{C}}(\ker (D_1)_3) \\ +\dim _{\mathbb{C}}(\ker 
(D_1')_3))-24(\dim\, G)-32(\dim\, H)-9(\dim \, L)=38.
\end{eqnarray*}
On the other hand $\dim\, \mathfrak{U}\ge 10\chi -2K^2=38$ by general principles (see
[Cat1b], p. 484). Thus $\dim\, \mathfrak{U}=38$ and one recovers the following theorem
(cf. [Ro\ss ], p. 116, thm. 1):
\begin{thm}
Regular surfaces of general type with $p_g=5,\: K^2=11$ such that the canonical
map is a birational morphism and the image $Y\subset \mathbb{P}^4$ has only improper
double points as singularities form an irreducible unirational open set $\mathfrak{U}$
of dimension $38$ inside their moduli space.
\end{thm}

\section[General $K^2$]{The case of general $K^2$}
\setcounter{equation}{0}
Recall that for a canonical surface $\pi:S\to Y\subseteq \mathbb{P}^4$
with
$p_g=5,\;q=0$ one has a resolution of the canonical ring
\begin{eqnarray}
\begin{CD}\mathbf{R}_{\bullet}:\;0@>>>\mathcal{A}(-6)\oplus\mathcal{A}(-4)^n@
>{{-\beta^t}\choose{\alpha^t}}>>\mathcal{A}(-3)^{2n+2}\end{CD}\nonumber\\
\begin{CD}@>{(\alpha\,\beta)}>>\mathcal{A}
\oplus\mathcal{A}(-2)^n@>>>\mathcal{R}@>>>0
\end{CD}
\end{eqnarray} where $n:=K^2-9$. In this section, as sort of a converse to
this, I propose to prove
\begin{thm}
Let $A=(\alpha\,\beta)$ be an $(n+1)\times (2n+2)$ matrix with entries in
the $1^{st}$ row cubic forms and linear forms otherwise such that
$\alpha\beta^t=\beta\alpha^t$ is symmetric. Put  $A':=(A$ with $1^{st}$
row erased) and
\begin{eqnarray*}M^s(n,2n+2):=\{[\mathbf{M}]:\mathbf{M}=(\mathbf{a}\,\mathbf{b})\in\mathbb{
C}^{n\times (2n+2)}\;
s.t.\;\mathbf{a}\mathbf{b}^t=\mathbf{b}\mathbf{a}^t\}\\ \subset\mathbb{P}^{n(2n+2)-1},
\end{eqnarray*}(where $\mathbf{a},\mathbf{b}\in\mathbb{C}^{n\times(n+1)}$) and \[
\Delta:=\{[\mathbf{M}]\in
M^s(n,2n+2):\mathrm{rk}\,\mathbf{M}\le n-1\}\subset M^s(n,2n+2).
\]Then the degeneracy locus $\Delta$ sits as an irreducible
4-codimensional subvariety inside the irreducible variety $M^s(n,2n+2)$.
Assume that the linear forms in $A'$ parametrically define a
$\mathbb{P}^4$ inside $M^s(n,2n+2)$ that is transverse to the locus
$\Delta$, the $(n+1)\times(n+1)$ minors of $A$ have no common factor and
that $\mathrm{Ann}_{\mathcal{A}}\,\mathcal{R}$ is prime,
where
$\mathcal{R}:=\mathrm{coker}(A:\mathcal{A}(-3)^{2n+2}\to\mathcal{A}\oplus\mathcal{A}(-2)^n)$.
Then $Y=\mathrm{Supp}\,\mathcal{R}\subset\mathbb{P}^4$, with closed subscheme
structure given by $\mathrm{Ann}_{\mathcal{A}}\,\mathcal{R}\subset\mathcal{A}$, is an
irreducible surface  with isolated nonnormal locus defined by
$I_n(A')$ as a reduced subscheme; assume further that the latter points are improper
double points of $Y$. Then
$\mathcal{R}$ gets a ring structure via
$\mathcal{R}=\mathrm{Hom}_{\mathcal{A}_Y}(\overline{I_n(A')}\cdot
\mathcal{A}_Y,\overline{I_n(A')}\cdot\mathcal{A}_Y)
$, where the bar denotes saturation.\\Moreover,
$X=\mathrm{Proj}(\mathcal{R})$ is then the canonical model of a surface
of general type with $K^2=n+9,\;p_g=5,\;q=0$ provided $X$ has only
rational double points as singularities. The morphism $\psi : X\to Y$ induced by the
inclusion $\mathcal{A}_Y\subset\mathcal{R}$ is finite and birational, and
$Y\subset\mathbb{P}^4$ is a canonical surface.
\end{thm}
\begin{proof}
To begin with, let's check that $M^s(n,2n+2)$ resp. $\Delta$ are
irreducible and $\mathrm{codim}_{M^s(n,,2n+2)}\Delta=4$ as stated. In
fact, I will prove that the irreducible algebraic group
$PGl_n(\mathbb{C})\times PSp_{2n+2}(\mathbb{C})$ acts (morphically) on
$M^s(n,2n+2)$ (via left resp. right multiplication) with orbits\\
$M_k^s(n,2n+2)-M_{k-1}^s(n,2n+2)$, $k=0,\ldots,n$, where
$M_k^s(n,2n+2)$ is the locus of matrices inside $M^s(n,2n+2)$ of rank
$\le k$. Since $M^s(n,2n+2)-M_{n-1}^s(n,2n+2)$ resp.
$M^s_{n-1}(n,2n+2)-M^s_{n-2}(n,2n+2)$ are clearly dense in
$M^s(n,2n+2)$ resp. $\Delta$, the latter are then irreducible.\\
Given $(\mathbf{a},\mathbf{b})\in M^s(n,2n+2)$ I will now use the
operations (Op) of Lemma 2.1 which belong to $PGl_n(\mathbb{C})\times
PSp_{2n+2}(\mathbb{C})$. Applying  (Op),(i) and (iv), possibly (v) one
transforms $(\mathbf{a},\mathbf{b})$ to get
$\left( \begin{array}{*{2}{c@{\:}}|*{1}{c@{\:}}} Id_r &0&  \\ 0&0&
\raisebox{1.5ex}[-1.5ex]{$\mathbf{b}'$}
\end{array}\right)$, where $Id_r$ is the $r\times r$ identity matrix, $r$
the rank of $\mathbf{a}$ and $\mathbf{b}'$ an $n\times(n+1)$ matrix.
Using the symmetry, one finds that one has actually put
$(\mathbf{a},\mathbf{b})$ in the shape
\[\left( \begin{array}{*{2}{c@{\:}}|*{2}{c@{\:}}} Id_r &0& \mathbf{b}_1'
& \mathbf{b}_2'
\\ 0&0&0 & \mathbf{b}_3'
\end{array}\right), \]where $\mathbf{b}_1',\:\mathbf{b}_2',\:\mathbf{b}_3'
$ are $r\times r,\:r\times (n+1-r),\:(n-r)\times(n+1-r)$ matrices,
respectively, and $\mathbf{b}_1'$ is symmetric. Using (Op), (iii) I get
\[\left( \begin{array}{*{2}{c@{\:}}|*{2}{c@{\:}}} Id_r &0& \mathbf{b}_1'
& 0
\\ 0&0&0 & \mathbf{b}_3'
\end{array}\right), \]and using (Op), (iii) again and the fact that
$\mathbf{b}_1'$ is symmetric: \[\left( \begin{array}{*{2}{c@{\:}}|*{2}{c@{\:}}} Id_r
&0& 0 & 0
\\ 0&0&0 & \mathbf{b}_3'
\end{array}\right). \]Finally, using (Op),(i),(iv) and (v), afterwards
(vi), one arrives at\\ $\left( \begin{array}{*{2}{c@{\:}}|*{1}{c@{\:}}}
Id_{r+s} &0&  \\ 0&0&
\raisebox{1.5ex}[-1.5ex]{0}
\end{array}\right)$, where $s$ is the rank of $\mathbf{b}_3'$. Thus the
orbits of the action of $PGl_n(\mathbb{C})\times PSp_{2n+2}(\mathbb{C})$
on $M^s(n,2n+2)$ are the ones mentioned above.\\
\\

Denote homogeneous coordinates in $\mathbb{P}^{n(2n+2)-1}$ by
$\{\mathbf{a}_{ij};\mathbf{b}_{ij}\}_{1\le i\le n,\\1 \le j\le n+1}$ and
consider the linear subspace \[ \Lambda:=\{\mathbf{a}_{ij}=0,1\le i\le
n,1\le j\le n+1;\:\mathbf{b}_{kl}=0,1\le k\le l\le n+1\}
\] and the restriction of the projection with center $\Lambda$
$\pi_{\Lambda}:M^s(n,2n+2)-\Lambda\to\mathbb{P}^N$, where
$N=n(n+1)+(n+1)(n+2)/2-2$. Clearly, $\pi_{\Lambda}$ is dominant and
generically one-to-one ($\pi_{\Lambda}|(\pi^{-1}_{\Lambda}U)$, where $U=\{ \det
(\mathbf{a}_{ij})_{1\le i,j\le h} \neq 0,\: h=1,\ldots ,n \}$, is one-to-one),
 hence
$\mathrm{dim}M^s(n,2n+2)=N$. Consider the following incidence correspondence with the
indicated two projections:\\
\setlength{\unitlength}{1cm}
\small
\begin{picture}(10,4)
\put(0,3){$\tilde{M}^s_{n-1}(n,2n+2):=\{([\mathbf{M}],L)\in
M^s(n,2n+2)\times Grass(n-1,n):\mathrm{im}(\mathbf{M})\subset L\}$}
\put(1,1){$M^s(n,2n+2)$}\put(9,1){$Grass(n-1,n)$}
\put(4.5,2.8){\vector(-2,-1){2.3}} \put(7,2.8){\vector(2,-1){2.5}}
\put(3.5,2){$pr_1$} \put(7.5,2){$pr_2$}
\end{picture}\normalsize
\\Then $pr_2$ is surjective, and choosing a suitable basis in
$\mathbb{C}^n$ one can identify the fibre of $pr_2$ over a point in
$Grass(n-1,n)$ with (the projectivisation of) the set of matrices of the
form $\left( \begin{array}{c|c} \tilde{\mathbf{a}} &
\tilde{\mathbf{b}}\\ 0 & 0\end{array}\right)$, where now
$\tilde{\mathbf{a}},\:\tilde{\mathbf{b}}$ are $(n-1)\times(n+1)$
matrices with $\tilde{\mathbf{a}}\tilde{\mathbf{b}}^t$ symmetric.
Analogously to the proof of the irreducibility of and computation of
the dimension of
$M^s(n,2n+2)$, one therefore finds that the fibres
of
$pr_2$ are irreducible and their dimension equals
$(n-1)(n+1)+(n+1)(n+2)/2-4$. Therefore
$\mathrm{dim}\tilde{M}^s_{n-1}(n,2n+2)=(n-1)(n+2)+(n+1)(n+2)/2-4$. Since
$pr_1$ is generically one-to-one onto
$\Delta$, one has
$\mathrm{codim}_{M^s(n,2n+2)}(\Delta)=n(n+1)+(n+1)(n+2)/2-2-(n-1)(n+2)-(n+1)(
n+2)/2+4=4.$\\
Now let $A=(\alpha\;\beta)$ be a matrix of forms meeting the
requirements of the theorem. Since the linear forms in $A'$ are
supposed to define a $\mathbb{P}^4$ inside $M^s(n,2n+2)$ transverse to
the locus $\Delta$, $I_n(A')$ (scheme-theoretically) defines a finite
set of reduced points in $\mathbb{P}^4$. I then have the fundamental
\begin{lemm}
If $A=(\alpha\;\beta)$ is an $(n+1)\times(2n+2)$ matrix with first row
cubic forms, other rows linear forms on $\mathbb{P}^4$ with
$\alpha\beta^t=\beta\alpha^t$ and such that $I_n(A')$ defines a set of
reduced points, then \begin{equation}
\overline{I_n(A')}=\overline{I_n(A)},\end{equation}where the bar denotes
saturation. 
\end{lemm}
\begin{proof}

Let $P$ be one of the points that $I_n(A')$ defines. I work locally, in
the ring of germs of regular functions around $P$,
$\mathcal{O}_{\mathbb{P}^4,P}$. Slightly abusing notation, I write again
$A$ for the matrix of the map
$\mathcal{O}^{2n+2}_{\mathbb{P}^4,P}\to\mathcal{O}^{n+1}_{\mathbb{P}^4,P}$
induced by $A$. I write
$A=\left(\begin{array}{c}a_1\\a_2\\A''\end{array}\right)$, where $a_1$
denotes the first row of $A$, $a_2$ the second, and $A''$ the
$(n-1)\times(2n+2)$ residual matrix. Replacing $a_2$ by a suitable linear
combination of the rows in $A'$, I may assume that the vanishing of the
entries in $a_2$ defines $P$. Since $P$ is reduced, then
$\mathrm{rank}A''=n-1$. At this point I again use the operations (Op) of
lemma 2.1 to transform $A''$ resp. $A$ such that I can easily read off
the information I want. During the calculation I will treat $A''$ resp.
$A$ as dynamical variables. I can find a unit among the entries of the
first row of $A''$ and using (Op), (vi) and (v) I can assume $A''_{11}$
is a unit. Using (Op), (iv) I reach
\[A''=\left( \begin{array}{*{4}{c@{\:}}|*{4}{c}}1&0&\ldots&0& & &\ast& \\
\hline  & & & & & & & \\ & & \raisebox{1.5ex}[-1.5ex]{$\ast$}& & &
&\raisebox{1.5ex}[-1.5ex]{$\ast$} & 
\end{array}\right)
\]and using (Op), (ii) I can assume $A''_{1,n+2}$ is a unit and using
(iv) I reach  $A''=\left(
\begin{array}{*{4}{c@{\:}}|*{4}{c@{\:}}}1&0&\ldots&0&u_1 &0 &\ldots&0 \\
\hline  & & & & & & & \\ & & \raisebox{1.5ex}[-1.5ex]{$\ast$}& & &
&\raisebox{1.5ex}[-1.5ex]{$\ast$} & 
\end{array}\right) $, where $u_1$ is a unit. Remarking that by symmetry,
$A''_{i,n+2}=u_1A''_{i,1},\;i=2,\ldots,n+1$, I can reach using row
operations
\[A''=\left(
\begin{array}{*{4}{c}|*{4}{c}}1&0&\ldots&0&u_1 &0 &\ldots&0 \\
 0& & & &0 & & & \\ \vdots& &\ast &
&\vdots& &\ast &
\\ 0& & & &0 & & & 
\end{array}\right), \]thus inductively
\[A=\left( \begin{array}{*{5}{c@{\:}}|*{5}{c@{\:}}}A_{11} & \cdots &
A_{1,n-1} & A_{1,n} & A_{1,n+1} & A_{1,n+2} & \cdots & A_{1,2n} &
A_{1,2n+1} & A_{1,2n+2}\\A_{21} & \cdots &
A_{2,n-1} & A_{2,n} & A_{2,n+1} & A_{2,n+2} & \cdots & A_{2,2n} &
A_{2,2n+1} & A_{2,2n+2}\\ \hline 1 & \cdots & 0 & 0 & 0 & u_1 & \cdots &
0 & 0 & 0\\ \vdots & \ddots & \vdots & \vdots & \vdots & \vdots & \ddots
& \vdots & \vdots & \vdots \\ 0 & \cdots & 1 & 0 & 0 & 0 & \cdots & u_n &
0 & 0
\end{array}\right).\]
Symmetry gives $A_{i,j}u_i=A_{i,n+1+j}$, $i=1,2\;j=1,\ldots,n-1$;
therefore, using (Op), (ii), I get
\[A=\left( \begin{array}{*{5}{c@{\:}}|*{5}{c@{\:}}}A_{11} & \cdots &
A_{1,n-1} & A_{1,n} & A_{1,n+1} & 0 & \cdots & 0 &
A_{1,2n+1} & A_{1,2n+2}\\A_{21} & \cdots &
A_{2,n-1} & A_{2,n} & A_{2,n+1} & 0 & \cdots & 0 &
A_{2,2n+1} & A_{2,2n+2}\\ \hline 1 & \cdots & 0 & 0 & 0 & 0 & \cdots &
0 & 0 & 0\\ \vdots & \ddots & \vdots & \vdots & \vdots & \vdots & \ddots
& \vdots & \vdots & \vdots \\ 0 & \cdots & 1 & 0 & 0 & 0 & \cdots & 0 &
0 & 0
\end{array}\right).\]
To ease notation, I put
$A_{1,n}=:A_1,\:A_{1,n+1}=:A_2,\:A_{1,2n+1}=A_3,\:A_{1,2n+2}=:A_4$ and
$A_{2,n}=:a_3,\:A_{2,n+1}=:a_4,\:A_{2,2n+1}=:-a_1,\:A_{2,2n+2}=:-a_2$ and
using row operations I get

\[
A=\left(
\begin{array}{*{4}{c@{\:}}|*{4}{c@{\:}}}
 & & A_1& A_2 & & & A_3 & A_4\\
\multicolumn{2}{c}{\raisebox{1.5ex}[-1.5ex]{0}} & a_3 & a_4 &
\multicolumn{2}{c}{\: \raisebox{1.5ex}[-1.5ex]{0}} & -a_1 & -a_2\\
 & & & & & & & \\ \multicolumn{2}{c}{\raisebox{1.5ex}[-1.5ex]{$I_{n-1}$}}
& \multicolumn{2}{c}{\raisebox{1.5ex}[-1.5ex]{0}} &
\multicolumn{2}{c}{\: \raisebox{1.5ex}[-1.5ex]{0}} &
\multicolumn{2}{c}{\raisebox{1.5ex}[-1.5ex]{0}}
\end{array}
\right)
\] and to make $A$ look more symmetric I can use (Op), (ii) to arrive at
\begin{equation}
A=\left(
\begin{array}{*{4}{c@{\:}}|*{4}{c@{\:}}}
 & & A_1& A_2 & & & A_3 & A_4\\
\multicolumn{2}{c}{\raisebox{1.5ex}[-1.5ex]{0}} & a_3 & a_4 &
\multicolumn{2}{c}{\: \raisebox{1.5ex}[-1.5ex]{0}} & -a_1 & -a_2\\
 & & & & & & & \\ \multicolumn{2}{c}{\raisebox{1.5ex}[-1.5ex]{$I_{n-1}$}}
& \multicolumn{2}{c}{\raisebox{1.5ex}[-1.5ex]{0}} &
\multicolumn{2}{c}{\: \raisebox{1.5ex}[-1.5ex]{$I_{n-1}$}} &
\multicolumn{2}{c}{\raisebox{1.5ex}[-1.5ex]{0}}
\end{array}
\right).
\end{equation}
Therefore near $P$, $A$ looks like (3); $a_1,\:a_2,\:a_3,\:a_4$ is a
regular sequence in $\mathcal{O}_{\mathbb{P}^4,P}$ since the second row
of $A$ was assumed to define $P$ and the property is invariant under the
above process. Moreover, $I_n(A')=\langle a_1,a_2,a_3,a_4 \rangle$ and
$I_n(A)=\langle A_1,A_2,A_3,A_4,\:a_1,a_2,a_3,a_4 \rangle$ near $P$. But
the symmetry $\sum\limits^{4}_{i=1} A_ia_i=0$ is a Koszul relation saying
$\exists$ a skewsymmetric matrix of elements $\{Q_{ij}\}_{1\le i,j\le 4}$
such that $A_i=\sum\limits_{j=0}^{4}Q_{ij}a_j,\;i=1,\ldots,4$. Therefore
$I_n(A')=I_n(A)$ locally near $P$.\\
Let me now again write $A$ resp. $A'$ for the given matrices of forms on
$\mathbb{P}^4$ globally. From the above, I find the equality of sheaves
$\widetilde{I_n(A')}=\widetilde{I_n(A)}$ (the latter symbols denoting the sheaves
associated to $I_n(A')$ resp. $I_n(A)$ on $\mathbb{P}^4$). Translated
back into the language of ideals this just says \begin{math}
\overline{I_n(A')}=\overline{I_n(A)}\end{math}.
\end{proof}
Now given the existence of a matrix $A=(\alpha \:\beta)$ with the
desirable properties as in the statement of the theorem I have a complex
\begin{eqnarray*}
\begin{CD}\mathbf{R}_{\bullet}:\;0@>>>\mathcal{A}(-6)\oplus\mathcal{A}(-4)^n@
>{{-\beta^t}\choose{\alpha^t}}>>\mathcal{A}(-3)^{2n+2}\end{CD}\\
\begin{CD}@>{(\alpha\,\beta)}>>\mathcal{A}
\oplus\mathcal{A}(-2)^n@>>>\mathcal{R}@>>>0
\end{CD}
\end{eqnarray*}and the requirement that the $(n+1)\times(n+1)$ minors of $A$ have no
common factor translates as
$\mathrm{grade}I_{n+1}(A)=\mathrm{codim}_{\mathcal{A}}I_{n+1}(A)\ge 2$
whence this complex is exact by the Eisenbud-Buchsbaum acyclicity
criterion (see e.g. [Ei], thm. 20.9, p. 500); moreover, since
$\sqrt{I_{n+1}(A)}= \sqrt{\mathrm{ann}_{\mathcal{A}}(\mathcal{R})}$ and then
$2\le\mathrm{grade}(I_{n+1}(A),\mathcal{A})=\mathrm{grade}(\mathrm{ann}_{\mathcal{A}}(\mathcal{R}),\mathcal{A})\equiv\mathrm{grade}
(\mathcal{R})\le \mathrm{projdim}_{\mathcal{A}}(\mathcal{R})=2$ (for the latter
inequality "projective dimension bounds grade" see e.g. [B-He], p. 25), I have 
$\mathrm{codim}_{\mathcal{A}}I_{n+1}(A)=2$ and
$Y:=\mathrm{Supp}(\mathrm{Proj}(\mathcal{A}/I_{n+1}(A)))=\mathrm{Supp}(\mathcal{R})$
is a surface in $\mathbb{P}^4$, irreducible by assumption.\\
I now intend to exploit lemma 3.2 to investigate the ring structure of
$\mathcal{R}$. First, quite generally I have \[
\mathcal{R}\subseteq \mathrm{Hom}_{\mathcal{A}_Y}(I_n(A')\cdot
\mathcal{A}_Y,I_n(A)\cdot \mathcal{A}_Y).
\]

For let me write $\{ 1,\:v_1,\ldots,v_n\}$ for the minimal set of generators of 
$\mathcal{R}$ corresponding to the standard basis of
$\mathcal{A}\oplus\mathcal{A}(-2)^n$ and let $M=$\\ $=(m_{ij})_{1\le i\le n+1,1\le
j\le n}$ be an arbitrary $(n+1)\times n$ submatrix of $A$, $M={m\choose M'}$,
where $m$ is the first row of $M$,
$M'=(m'_{ij})_{1\le i,j\le n}$ the $n\times n$ residual matrix. Let $v_k$ be
one of the $\{v_1,\ldots v_n\}$. In $\mathcal{R}$ I have relations\[
m_{1i}+\sum\limits_{j=1}^nm'_{ji}v_j=0,\;i=1,\ldots,n
\]hence also\[
(m_{1i}+\sum\limits_{j=1}^n m'_{ji}v_j)M^{'\ast}_{ik}=0,\;i=1,\ldots,n
\]denoting by $M^{'\ast}_{ik}$ the $i,k-$entry of the adjoint $M^{'\ast}$ of
$M'$. Adding the latter equations up for the various $i$ gives
$\det(M')v_k=\pm\det(M$  with $k^{th}$ row deleted$)$. This shows that the
stated inclusion holds (basically as a consequence of Cramer's rule). By lemma
3.2 it follows that $\mathcal{R}\subseteq \mathrm{Hom}_{\mathcal{A}_Y}(I_n(A')\cdot
\mathcal{A}_Y,\overline{I_n(A')}\cdot \mathcal{A}_Y).$\\
I claim that then also  $\mathcal{R}\subseteq
\mathrm{Hom}_{\mathcal{A}_Y}(\overline{I_n(A')}\cdot
\mathcal{A}_Y,\overline{I_n(A')}\cdot \mathcal{A}_Y).$ In fact,
$\overline{I_n(A')}=\{p\in\mathcal{A}=\mathbb{C}[x_0,\ldots,x_4] :$ for each
$i=0,\ldots,4\:\exists n$ such that $x_i^np\in I_n(A')\}$. But then for
$\bar{p}\in \overline{I_n(A')}\cdot \mathcal{A}_Y,\;v\in \mathcal{R},$ the
expression $\bar{p}v$ is again in $ \overline{I_n(A')}\cdot \mathcal{A}_Y$: For
$i$ among $0,\ldots,4$ $(x_i^n\bar{p})v=x_i^n(\bar{p}v)$ is in
$\overline{I_n(A')}\cdot
\mathcal{A}_Y$, therefore there exists an integer $m$ such that
$x_i^mx_i^n(\bar{p}v)\in I_n(A')\cdot \mathcal{A}_Y$, i.e.
$\bar{p}v\in\overline{I_n(A')}\cdot \mathcal{A}_Y$.\\
Therefore I get the chain of inclusions\begin{equation}
\mathcal{R}\subseteq
\mathrm{Hom}_{\mathcal{A}_Y}(\overline{I_n(A')}\cdot\mathcal{A}_Y,\overline{I_n(A')}\cdot
\mathcal{A}_Y)\subseteq
\mathrm{Hom}_{\mathcal{A}_Y}(\overline{I_n(A')}\cdot\mathcal{A}_Y,\mathcal{A}_Y).
\end{equation}
To show the reverse inclusion
$\mathrm{Hom}_{\mathcal{A}_Y}(\overline{I_n(A')}\cdot\mathcal{A}_Y,\mathcal{A}_Y)\subseteq
\mathcal{R}$ I need another technical result. Let me introduce the so called
conductor $\mathcal{C}$ of $\mathcal{R}$ into $\mathcal{A}_Y$,
$\mathcal{C}:=\mathrm{Hom}_{\mathcal{A}_Y}(\mathcal{R},\mathcal{A}_Y)$, and the
associated sheaf on $\mathbb{P}^4$,
$\tilde{\mathcal{C}}:=\mathcal{H}om_{\mathcal{O}_Y}(\tilde{\mathcal{R}},\mathcal{O}_Y)$. 
 \begin{lemm}
$\mathcal{R}$ being as in the statement of the theorem, one has
$\mathcal{R}=\Gamma_\ast(\tilde{\mathcal{R}})$, where $\tilde{\mathcal{R}}$ is the
sheaf on $\mathbb{P}^4$ associated to the graded module $\mathcal{R}$, supported on
the surface $Y$. Moreover the fact that locus defined by
$I_n(A')$ as a reduced subscheme is a finite number of points which are improper
double points on $Y$, implies
that
$\tilde{\mathcal{R}}$ is reflexive in the sense that
$\tilde{\mathcal{R}}=\mathcal{H}om_{\mathcal{O}_Y}(\tilde{\mathcal{C}},\mathcal{O}_Y)$.
\end{lemm}
\begin{proof}First, $\mathcal{R}$ equals the full module of sections of
$\tilde{\mathcal{R}}$, i.e.
$\mathcal{R}=\Gamma_\ast(\tilde{\mathcal{R}})$. For put
$\mathcal{F}_0:=\mathcal{O}\oplus\mathcal{O}(-2)^n$,
$\mathcal{F}_1:=\mathcal{O}(-3)^{2n+2}$,
$\mathcal{F}_2:=\mathcal{O}(-6)\oplus\mathcal{O}(-4)^n$ and sheafify the resolution
of $\mathcal{R}$ above to get the
diagram\\ 
\setlength{\unitlength}{1cm}
\begin{picture}(0,4)
\put(3,3){$\mathcal{F}_2$}\put(1.5,3){0} \put(4.5,3){$\mathcal{F}_1$}
\put(6.5,3){$\mathcal{F}_0$}\put(8,3){$\tilde{\mathcal{R}}$}\put(9.5,3){0}
\put(5.5,2){$\mathcal{G}$}
\put(1.8,3.12){\vector(1,0){0.9}}\put(3.4,3.12){\vector(1,0){0.9}}
\put(4.9,3.12){\vector(1,0){1.4}}\put(6.9,3.12){\vector(1,0){0.9}}
\put(8.4,3.12){\vector(1,0){0.9}}\put(4.85,2.8){\vector(1,-1){0.6}}
\put(5.75,2.25){\vector(1,1){0.6}}
\put(6.45,1.2){0}\put(4.8,1.2){0}\put(5.05,1.5){\vector(1,1){0.5}}
\put(5.8,1.9){\vector(1,-1){0.5}}
\end{picture}
From this one gets the exact sequences
\[0\rightarrow H^0(\mathbb{P}^4,\mathcal{G}(j))\rightarrow
H^0(\mathbb{P}^4,\mathcal{F}_0(j))\rightarrow
H^0(\mathbb{P}^4,\mathcal{R}(j))\rightarrow 0 \]
\[0\rightarrow H^0(\mathbb{P}^4,\mathcal{F}_2(j))\rightarrow
H^0(\mathbb{P}^4,\mathcal{F}_1(j))\rightarrow
H^0(\mathbb{P}^4,\mathcal{G}(j))\rightarrow 0 \]since
$\forall\;j\; H^1(\mathbb{P}^4,\mathcal{G}(j))=0$,
$H^1(\mathbb{P}^4,\mathcal{F}_2(j))=0$  ($H^i(\bigoplus\mathcal{O}(d_k))=0
$, $i\neq0,4 $ on $\mathbb{P}^4$). Putting the above two exact sequences
together  gives

\[\begin{CD}
0\rightarrow\Gamma_j(\mathcal{F}_2)@>>>\Gamma_j(\mathcal{F}_1)@>>>\Gamma_j(
\mathcal{F}_0)@>>>\Gamma_j(\tilde{\mathcal{R}})\rightarrow 0\\
@VV{\cong}V  @VV{\cong}V  @VV{\cong}V  @VV{\iota}V\\
0\rightarrow F_2(j) @>>> F_1(j) @>>> F_0(j) @>>> \mathcal{R}(j)
\rightarrow 0
\end{CD}\]
and $\iota$ is an isomorphism (sc. $F_2,\:F_1,\:F_0$ the graded free
modules appearing in the resolution (1)).\\
Secondly,
$\tilde{\mathcal{R}}=\mathcal{H}om_{\mathcal{O}_Y}(\tilde{\mathcal{C}},\mathcal{O}_Y)$,
where
$\tilde{\mathcal{C}}:=\mathcal{H}om_{\mathcal{O}_Y}(\tilde{\mathcal{R}},\mathcal{O}_Y)$
is the sheaf of conductors of $\tilde{\mathcal{R}}$ into $\mathcal{O}_Y$.
Namely, for $P$ a point where $I_n(A')$ does not drop rank the natural homomorphism
$\tilde{\mathcal{R}}_P\to\mathrm{Hom}_{\mathcal{O}_{Y,P}}(\tilde{\mathcal{C}}_P,\mathcal{O}_{Y,P})$
is clearly an isomorphism because then locally at $P$ $A'$ is surjective, and
$\mathcal{R}/\mathcal{A}_Y$ being the cokernel of the matrix $A'$, I have
$\tilde{\mathcal{R}}_P=\mathcal{O}_{Y,P}$ and also
$\tilde{\mathcal{C}}_P=\mathcal{O}_{Y,P}$. Therefore the interest is in the
improper double points of $Y$. \\
Therefore let $Q$ be one of the improper double points that $I_n(A')$ defines. Then
locally around $Q$ $A$ can be written as
\[
A=\left(
\begin{array}{*{4}{c@{\:}}|*{4}{c@{\:}}}
 & & A_1& A_2 & & & A_3 & A_4\\
\multicolumn{2}{c}{\raisebox{1.5ex}[-1.5ex]{0}} & a_3 & a_4 &
\multicolumn{2}{c}{\: \raisebox{1.5ex}[-1.5ex]{0}} & -a_1 & -a_2\\
 & & & & & & & \\ \multicolumn{2}{c}{\raisebox{1.5ex}[-1.5ex]{$I_{n-1}$}}
& \multicolumn{2}{c}{\raisebox{1.5ex}[-1.5ex]{0}} &
\multicolumn{2}{c}{\: \raisebox{1.5ex}[-1.5ex]{$I_{n-1}$}} &
\multicolumn{2}{c}{\raisebox{1.5ex}[-1.5ex]{0}}
\end{array}
\right)
\]
as was shown above (see (3)). I have that $\tilde{\mathcal{C}}_Q=(a_1,a_2,a_3,a_4)$
because $\varphi\in\mathcal{O}_{Y,Q}$ is in $\tilde{\mathcal{C}_Q}\:\Leftrightarrow
\:\exists p,q\in\mathcal{O}_{\mathbb{P}^4,Q}$ s.t. $pv+q=0$ in
$\tilde{\mathcal{R}}_Q$ and $p$ is a lift of $\varphi$ ($1,v$ denoting the minimal
set of generators in $\tilde{\mathcal{R}}_Q$ corresponding to the first two rows of
$A$ as above). Since
$a_1,\:a_2,\:a_3,\:a_4$ define the improper double point $Q$ as a reduced subscheme,
I can (without loss of generality) assume that $a_i=x_i,\: i=1,\ldots,4$ are
coordinates in
$\mathbb{C}[[x_1,\ldots ,x_4]]$ and
$\mathcal{O}_{Y,Q}^{an}=\mathbb{C}[[x_1,\ldots,x_4]]/(x_1,x_2)\cap (x_3,x_4)$,
changing to the analytic category; then
$\tilde{\mathcal{C}}_Q^{an}=(x_1,x_2,x_3,x_4)
\mathcal{O}_{Y,Q}=(x_1,x_2,x_3,x_4)/(x_1x_3,x_2x_3,x_1x_4,x_2x_4)=(x_1,x_2)\mathbb{C}
[[x_1,x_2]]\oplus (x_3,x_4)\mathbb{C}[[x_3,x_4]]$. I have to show
$\tilde{\mathcal{R}}_Q^{an}\subset$ \\ $\subset\mathrm{Hom}_{\mathcal{O}_{Y,Q}^{an}}
(\tilde{\mathcal{C}}_Q^{an},\mathcal{O}_{Y,Q}^{an})=\tilde{\mathcal{R}}_Q^{an\ast\ast}$
 is an isomorphism. Look at\\
\setlength{\unitlength}{1cm}
\begin{picture}(10,4)
\put(3,3){$\mathcal{O}_{Y,Q}^{an}$}
\put(3.9,3){$\stackrel{i}{\hookrightarrow}$}
\put(4.5,3){$\tilde{\mathcal{R}}_{Q}^{an}$}
\put(5.2,3){$\longrightarrow$}
\put(6.1,3){coker $i$}
\put(7.5,3){$\longrightarrow$}
\put(8.3,3){0}

\put(3.1,2.3){$\parallel$}
\put(4.6,2.3){$\bigcap$}
\put(6.4,2.3){$\bigcap$}

\put(3,1.6){$\mathcal{O}_{Y,Q}^{an}$}
\put(3.9,1.6){$\stackrel{i'}{\hookrightarrow}$}
\put(4.5,1.6){$\tilde{\mathcal{R}}_{Q}^{an\ast\ast}$}
\put(5.2,1.6){$\longrightarrow$}
\put(6.1,1.6){coker $i'$}
\put(7.5,1.6){$\longrightarrow$}
\put(8.3,1.6){0}

\put(4.6,1.2){$\parallel$}
\put(4,0.7){$\mathrm{Hom}_{\mathcal{O}_{Y,Q}^{an}}
(\tilde{\mathcal{C}}_Q^{an},\mathcal{O}_{Y,Q}^{an}).$}
\end{picture}\\
Here coker $i\neq 0$ since $\tilde{\mathcal{R}}_{Q}^{an}$ is minimally generated by 2
elements as an $\mathcal{O}_{Y,Q}^{an}-$module ( $\tilde{\mathcal{R}}_{Q}$ is  the
cokernel of the matrix in (3)). But on the other hand, one computes $\mathrm{Hom}_{\mathcal{O}_{Y,Q}^{an}}
(\tilde{\mathcal{C}}_Q^{an},\mathcal{O}_{Y,Q}^{an})=\mathbb{C}[[x_1,x_2]]\oplus
\mathbb{C}[[x_3,x_4]]$. Namely, $\mathbb{C}[[x_1,x_2]]\oplus
\mathbb{C}[[x_3,x_4]]\subset \mathrm{Hom}_{\mathcal{O}_{Y,Q}^{an}}
(\tilde{\mathcal{C}}_Q^{an},\mathcal{O}_{Y,Q}^{an})$ via
$(\varphi_1\oplus\varphi_2)(c_1\oplus c_2)=\varphi_1 c_1+\varphi_2 c_2$ for
$(\varphi_1\oplus\varphi_2)\in \mathbb{C}[[x_1,x_2]]\oplus
\mathbb{C}[[x_3,x_4]]$ and $(c_1\oplus c_2)\in (x_1,x_2)\mathbb{C}
[[x_1,x_2]]\oplus (x_3,x_4)\mathbb{C}[[x_3,x_4]]$, and conversely, given $\varphi\in \mathrm{Hom}_{\mathcal{O}_{Y,Q}^{an}}
(\tilde{\mathcal{C}}_Q^{an},\mathcal{O}_{Y,Q}^{an})$,
then \\$\varphi (x_{1/2})\in
\mathrm{Ann}(x_3)\cap\mathrm{Ann}(x_4)=\mathbb{C}[[x_1,x_2]]$ and $\varphi
(x_1)x_2=\varphi (x_2)x_1 \in \mathbb{C}[[x_1,x_2]]$ whence $\varphi
(x_{1/2})=\varphi ' x_{1/2}$ with $\varphi ' \in \mathbb{C}[[x_1,x_2]]$. Similarly,
$\varphi (x_{3/4})=\varphi '' x_{3/4}$ with $\varphi '' \in \mathbb{C}[[x_3,x_4]]$. Now
from the exact sequence
\[
0\to \mathbb{C}[[x_1,\ldots,x_4]]/(x_1,x_2)\cap (x_3,x_4)\to \mathbb{C}[[x_1,x_2]]\oplus
\mathbb{C}[[x_3,x_4]]\to \mathbb{C}\to 0
\]
I have coker $i'\cong \mathbb{C}$ whence the righthand inclusion in the above diagram
is an isomorphism and therefore also the middle one (by the 5-lemma).
\end{proof}
Now clearly $(I_n(A')\cdot\mathcal{A}_Y)^{\sim}\cong\tilde{\mathcal{C}}$ as
sheaves; for $I_n(A')\subseteq\mathrm{ann}_{\mathcal{A}}(\mathcal{R}/\mathcal{A}_Y)$
(Fitting's lemma), hence
$I_n(A')\cdot\mathcal{A}_Y\subseteq\mathcal{C}$ which gives me a morphism of these
sheaves which is an isomorphism (in fact, in the proof of the preceding lemma I saw
that for $P$ one of the points where
$I_n(A')$ drops rank $(I_n(A')\cdot\mathcal{A}_Y)^{\sim}_P$ and
$\tilde{\mathcal{C}}_P$ are both $(a_1,a_2,a_3,a_4)\mathcal{O}_{Y,P}$ and are
$\mathcal{O}_{Y,P}$ otherwise). Hence also
$\tilde{\mathcal{R}}=\mathcal{H}om_{\mathcal{O}_Y}((I_n(A')\cdot\mathcal{A}_Y)^{\sim},\mathcal{O}_Y)$,
 but the full module of sections of the latter sheaf contains
$\mathrm{Hom}_{\mathcal{A}_Y}(\overline{I_n(A')}\cdot\mathcal{A}_Y,\mathcal{A}_Y)$
and combined with the fact that $\tilde{\mathcal{R}}$ equals its full module of
sections and equation (4), I arrive at the fact that $\mathcal{R}$ is a ring via $\mathcal{R}=\mathrm{Hom}_{\mathcal{A}_Y}(\overline{I_n(A')}\cdot
\mathcal{A}_Y,\overline{I_n(A')}\cdot\mathcal{A}_Y)
$. Then $X=\mathrm{Proj}(\mathcal{R})$ is the canonical model of a surface of
general type $S$ with $q=0,\:p_g=5,\:K^2=n+9$ by theorem 1.6, if $X$ has only
rational double points as singularities.\\
The fact in the statement of theorem 3.1 that
$I_n(A')$ gives precisely the nonnormal locus of the surface $Y$ now follows from
remark 2 after the proof of theorem 1.6. This completes the proof of theorem 3.1.  
\end{proof}

\bfseries{Remark 1.} \mdseries
The variety $M^s(n,2n+2)$ is a complete intersection of $\frac{n(n-1)}{2}$ quadrics;
namely by the argument at the beginning  of the proof of theorem 3.1,
$\mathrm{codim}_{\mathbb{P}^{n(2n+2)-1}}(M^s(n,2n+2))=\frac{n(n-1)}{2}$ and
$M^s(n,2n+2)$ is cut out by the ${n \choose 2}=\frac{n(n-1)}{2}$ quadratic equations
given by the symmetry $\mathbf{a}\mathbf{b}^t=\mathbf{b}\mathbf{a}^t$. Moreover,
$\mathbb{P}^4$'s clearly exist on this variety. In fact, $\dim
Grass(5,n(2n+2))=5(n(2n+2)-5)$, and the Fano variety of 4-planes lying on one of the
above quadrics has codimension 15 in $Grass(5,n(2n+2))$ whence one finds an at least
$5(n(2n+2)-5)-15\frac{n(n-1)}{2}=\frac{5}{2}[n^2+7n-10]$ dimensional family of
$\mathbb{P}^4$'s on $M^s(n,2n+2)$. For $n\ge 2$ ($\Leftrightarrow K^2\ge 11$) this
number is positive. However, this says of course nothing as to whether one can find
$\mathbb{P}^4$'s transverse to the locus $\Delta$ of theorem 3.1. in all cases. Thus
the next step towards the construction of canonical surfaces with higher $K^2$,
$K^2\ge 13$, say, should be a more detailed analysis of the Fano variety of
$\mathbb{P}^4$'s on $M^s(n,2n+2)$.

\bfseries{Remark 2.}\mdseries
As a second step towards understanding the afore-mentioned surfaces one can look at
the forgetful maps:
\begin{eqnarray*}
F_n^{(1)}: \left\{
\begin{array}{l}
\mathrm{Parameter\: space\: of\:}(n+1)\times (2n+2) \mathrm{\:matrices\:}\\
(\alpha\:\beta) \mathrm{\:with\: the \: properties\: listed \: in}\\
\mathrm{the \: hypotheses\: of \: theorem\: 3.1}
\end{array}
\right\} \\ \longrightarrow 
\left\{
\begin{array}{l}
\mathrm{Parameter\: space\: of\:} n\times (2n+2) \mathrm{\: matrices\:}\\
(\alpha ' \: \beta ') \mathrm{\: of \: linear \: forms\: on \: }\mathbb{P}^4
\mathrm{ \: with\:}
 \alpha '
\beta^{'t} \\ 
\mathrm{\: symmetric\: and \: with \:} (\alpha ' \: \beta ') \mathrm{\: degenerating} 
\\
\mathrm{in \: a\: finite \: number \: of \: reduced\: points\: in\:} \mathbb{P}^4 
\end{array}
\right\}
\end{eqnarray*}
obtained by erasing the first row of a matrix $( \alpha \: \beta)$, and try to
understand 1) when $F^{(1)}_n$ is dominant, 2) what its fibres look like. Forgetting
even more, one can ask the same questions for the maps
\begin{eqnarray*}
F^{(2)}_n: \left\{
\begin{array}{l}
\mathrm{Parameter\: space\: of\:}(n+1)\times (2n+2) \mathrm{\:matrices\:}\\
(\alpha\:\beta) \mathrm{\:with\: the \: properties\: listed \: in}\\
\mathrm{the \: hypotheses\: of \: theorem\: 3.1}
\end{array}
\right\} \\ \longrightarrow \bigsqcup\limits_{\mathrm{const.\: polynomials\:}P}
Hilb^P_{\mathbb{P}^4}
\end{eqnarray*}
where the latter denotes the Hilbert scheme of points in $\mathbb{P}^4$, and the map
is given by sending a matrix $A=(\alpha\:\beta)$ to the points $I_n(A')$ defines.

\bfseries{Remark 3.}\mdseries
Finally, it would be interesting to find a purely algebraic proof of lemma 3.3, thus
going beyond the assumption that the nonnormal locus of $Y$ consists of improper
double points alone; namely, for $Q$ one of the points that $I_n(A')$ defines, I
have to prove that the $\mathcal{O}_{Y,Q}-$module $\tilde{\mathcal{R}}_Q$ is
reflexive. A possible strategy to see this algebraically is as follows: Locally at
$Q$, $A$ can be written as in (3). Put $B_1:=a_3,\:B_2:=a_4,\:B_3:=-a_1,B_4:=-a_2$
and note that $\mathcal{O}_{Y,Q}=\mathcal{O}_{\mathbb{P}^4,Q}/(A_iB_j-A_jB_i)_{1\le
i,j \le 4}$. For clearly $(A_iB_j-A_jB_i)\subset\mathrm{ann}(\tilde{\mathcal{R}}_Q)$
by Fitting's lemma, and conversely, writing $\{ 1,v\}$ for the minimal set of
generators of $\tilde{\mathcal{R}}_Q$, if $R\cdot 1=0$ in $\tilde{\mathcal{R}}_Q$,
$R\in \mathcal{O}_{\mathbb{P}^4,Q}$, then, by the symmetry, $\exists \lambda_i,i=1,
\ldots ,4:\: \sum \lambda_i B_i=0,\:\sum \lambda_i A_i=R$. The former relation is a
Koszul relation saying $\exists \mu_{ij}=-\mu_{ji},\: i,j\in \{ 1,\ldots ,4 \}:\:
\lambda_i=\sum \mu_{ij}B_j$. Therefore $R=\sum\limits_{i<j}\mu_{ij}(B_jA_i-B_iA_j)$,
whence
$(A_iB_j-A_jB_i)=\mathrm{ann}_{\mathcal{O}_{\mathbb{P}^4,Q}}(\tilde{\mathcal{R}}_Q)$.

Next, having the inclusion $\tilde{\mathcal{R}}_Q\subset
\mathrm{Hom}_{\mathcal{O}_{Y,Q}}(\tilde{C}_Q,\mathcal{O}_{Y,Q})$, I want to show that
every $\varphi\in \mathrm{Hom}_{\mathcal{O}_{Y,Q}}(\tilde{C}_Q,\mathcal{O}_{Y,Q})$
comes from an element in $\tilde{\mathcal{R}}_Q$ (recall
$\tilde{\mathcal{C}}_Q=(B_1,B_2,B_3,B_4)\mathcal{O}_{Y,Q}$ and
$B_iv=-A_i,\:B_i\cdot 1=B_i$). Then putting $\varphi (B_i)=:\beta_i$,
$A_i\beta_j=\varphi (A_iB_j)=\varphi (A_jB_i)=A_j\beta_i$ in $\mathcal{O}_{Y,Q}$.
Therefore I would get what I want if $(A_i\beta_j-A_j\beta_i)\subset (A_iB_j-A_jB_i)$
implied that the vector $(\beta_1,\ldots,\beta_4)$ is a linear combination (mod
$(A_iB_j-B_iA_j)$) of the vectors
$(A_1,\ldots, A_4)$ and $(B_1,\ldots,B_4)$; in other words, if the complex
\[
\begin{CD}
\mathcal{O}_{Y,Q}^2 @>{ \small \left( \begin{array}{*{2}{c@{\:}}} A_1 & B_1 \\
 A_2 & B_2 \\ A_3 & B_3  \\ A_4 & B_4 \end{array} 
\right)}>> \normalsize \mathcal{O}_{Y,Q}^4 @>{ \small \left(
\begin{array}{*{4}{c@{\:}}} A_2 &  -A_1 & 0 & 0
\\ A_3 & 0 & -A_1 & 0 \\ A_4 & 0 & 0 & -A_1 \\  0 & A_3 & -A_2 & 0 \\ 0 & A_4 & 0 &
-A_2 \\ 0 & 0 & A_4 & -A_3 
\end{array}
\right)}>>
\normalsize \mathcal{O}_{Y,Q}^6
\end{CD}
\]
was exact. If the $A_i,\: B_i$ are replaced with indeterminates $X_i,\: Y_i,\;
i=1,\ldots ,4$, over $\mathcal{O}_{\mathbb{P}^4,Q}$, then this is exact as I checked
using the computer algebra system MACAULAY2 (one should of course find a
theoretically satisfactory reason for this). Thus it would be nice to know exactly
which genericity assumptions on the $A_i,\: B_i$ are needed for the above
complex to remain exact when I specialize $X_i\leadsto A_i,\: Y_i\leadsto B_i$.

\section[A commutative algebra lemma]{A commutative algebra lemma}
\setcounter{equation}{0}
This section stands somewhat isolated from the rest of the treatise. I included it
merely to fix up a fact that slightly improves on a theorem of M. Grassi ([Gra]). To
find some amelioration of the structure theorem for Gorenstein algebras in
codimension 2 presented in [Gra] was actually the superordinate aim from which
this work departed.\\
I'd like to work in the generality and setting adopted in [Gra], so let:\\
$(R,\mathfrak{m},k):=$a Cohen-Macaulay local ring with $2\notin \mathfrak{m}$,\\
$A:=$a codimension 2 Gorenstein algebra over $R$, i.e. a finite $R-$algebra with 
$\mathrm{dim}(R)-\mathrm{dim}_R(A)=2$ and $A\cong \mathrm{Ext}_R^2(A,R)$ as
$A-$modules.\\
Finally, it will be convenient to have the concept of Koszul module available.
Whereas the usual Koszul complex is associated with a linear form $f: R^n\to R$, a
Koszul module is a module having a resolution similar to the Koszul complex up to
the fact that the r\^{o}le of $f$ is taken by a family of (vector-valued) maps from
$R^n$ to $R^n$. I'll only make this precise in the relevant special case:\\
A finite $R-$module $M$ having a length 2 resolution 
\begin{equation}
\begin{CD}
0 \to R^n @>{\rho_1\choose\rho_2}>> R^{2n} @>(\tau_1\:\tau_2)>> R^n \to M \to 0
\end{CD}
\end{equation}
some $n\in\mathbb{N}$, is a Koszul module iff
$\mathrm{det}(\tau_1),\;\mathrm{det}(\tau_2)$ is a regular sequence on $R$ and
$\exists$ a unit
$\lambda\in R$:
$\mathrm{det}(\rho_1)=(-1)^n\lambda
\mathrm{det}(\tau_2),\;\mathrm{det}(\rho_2)=\lambda \mathrm{det}(\tau_1)$.\\
Then Grassi proves in case $R$ is a domain ([Gra], thm. 3.3) that $A$ has a
(Gorenstein) symmetric resolution
\begin{equation}
\begin{CD}
0 \to R^n @>{-\beta^t\choose\alpha}>> R^{2n} @>(\alpha\:\beta)>> R^n \to A \to 0
\end{CD}
\end{equation}
and a second resolution of the prescribed type (1) for the Koszul module
condition, and that these 2 are related by an isomorphism of complexes which is
the identity in degrees 0 and 2; firstly, for sake of generality, I will briefly show
that the assumption "$R$ a domain" is in fact not needed, and secondly, prove that
there is one single resolution of $A$ meeting both requirements, i.e. a resolution
as in (2) with $\mathrm{det}(\alpha),\;\mathrm{det}(\beta)$ an $R-$regular
sequence. This still gives no indication of how the ring structure of $A$ is
encoded in the resolution, but as the concepts of Koszul module and Gorenstein
symmetric resolution seem to provide a pleasing setting to investigate this
question, it can be useful to have a result linking these two.\\
For the first part, I note that the only place in [Gra] where the hypothesis that
$R$ be a domain enters is at the beginning of the proof of proposition 1.5, page
930: Here one is given a resolution as in (1), but \underline{without} any
additional assumptions on
$\mathrm{det}(\tau_1),\mathrm{det}(\tau_2),\mathrm{det}(\rho_1),\mathrm{det}(\rho_2)$
whatsoever, and Grassi wants to conclude that $\exists$ a base change in $R^{2n}$
such that (in the new base) $\mathrm{det}(\tau_1)$ is not a zero divisor on $R$.
But this can be proven by a similar method as Grassi uses in the sequel of the
proof of proposition 1.5, without using "$R$ a domain":
For let $\mathfrak{p}_1,\ldots,\mathfrak{p}_r$ be the associated primes of $R$
which are precisely the minimal elements of $\mathrm{Spec}(R)$ since $R$ is CM.
One shows that $\exists$ a base change in $R^{2n}$ such that
$\mathrm{det}(\tau_1)\notin\mathfrak{p}_i,\;\forall i=1,\ldots,r$ (in the new
base), more precisely, that $\exists$ a sequence of $r$ base changes such that
after the $m$th base change
\[
(\ast) \hspace{1cm} \mathrm{det}(\tau_1)\notin\mathfrak{p}_i,\;\forall
i\in\{r-m+1,\ldots,r\},
\] 
$m=0,\ldots,r$, the assertion being empty for $m=0$. Therefore, inductively,
suppose $(\ast)$ holds for $m$ to get it for $m+1$.\\
Denote by $[i_1,\ldots,i_n]$ the maximal minor of $(\tau_1\;\tau_2)$
corresponding to the columns $i_1,\ldots,i_n,\;i_j\in\{1,\ldots,2n\}$. If
$[1,\ldots,n]\notin\mathfrak{p}_{r-m}$ I'm already O.K., so suppose
$[1,\ldots,n]\in\mathfrak{p}_{r-m}$. By the Eisenbud-Buchsbaum acyclicity
criterion $I_n((\tau_1\;\tau_2))$ cannot consist of zerodivisors on $R$ alone,
therefore set 
\begin{eqnarray*}
l_1:=\mathrm{min}\{c :\:\exists s_1,\ldots,s_{n-1}\:\mathrm{with}\:
s_1<s_2<\ldots <s_{n-1}<c \\ \mathrm{and}  
 [s_1,\ldots ,s_{n-1},c]\notin\mathfrak{p}_{r-m}\} \end{eqnarray*} (then $n<l_1\le
2n$) and inductively,
\begin{eqnarray*}
l_i:=\mathrm{min}\{c :\:\exists s_1',\ldots,s_{n-i}'\:\mathrm{with}\:
s_1'<\ldots<s_{n-i}'<c<l_{i-1}<\ldots<l_1\\
\mathrm{and}\:[s_1',\ldots,s_{n-i}',c,l_{i-1},\ldots,l_1]\notin\mathfrak{p}_{r-m}\}
,\end{eqnarray*} 
$i=2,\ldots,n$. Then $\exists \;J$ such that $n<l_J<l_{J-1}<\ldots<l_1\le 2n$ and
for $I>J\;l_I\in\{1,\ldots,n\}$ ($J=n$ might occur and then the set of
$l_I\in\{1,\ldots,n\}$ is empty; this does not matter). \\
I have $[l_n,\ldots,l_1]\notin\mathfrak{p}_{r-m}$ by construction. Choose
$b\in(\bigcap\limits_{i=r-m+1}^{r}\mathfrak{p}_i)\backslash\mathfrak{p}_{r-m}$,
which is nonempty since the $\mathfrak{p}_i$'s are the minimal elements of
$\mathrm{Spec}(R)$. Denote by $y_1< \ldots <y_J$ the complementary indices of the
$l_I\in\{1,\ldots ,n\}$ inside $\{ 1,\ldots ,n \}$ and consider the base change on
$R^{2n}$: $M_{y_1,l_J}(b)\circ M_{y_2,l_{j-1}}(b)\circ \ldots \circ
M_{y_J,l_1}(b)$, where $M_{y_{\nu},l_{J-\nu +1}}(b),\; \nu=1,\ldots ,J$ is addition
of $b$ times the $l_{J-\nu +1}$ column to the $y_{\nu}$ column. Then one sees (by
the multilinearity of determinants) 
\[
[1,\ldots,n]_{\mathrm{new}}=[1,\ldots,n]_{\mathrm{old}}\pm
b^J[l_n,\ldots,l_1]_{\mathrm{old}}+b\mu,
\]
where "new" means after and "old" before the base change and $\mu$ is an
element in $\mathfrak{p}_{r-m}$ by the defining minimality property of the $l$'s.
Therefore, since by the induction hypothesis $[1,\ldots,n]_{\mathrm{old}}\notin
\mathfrak{p}_i,\; \forall i\in\{r-m+1,\ldots,r\}$ and $b$ is chosen appropriately:
$[1,\ldots,n]_{\mathrm{new}}\notin \mathfrak{p}_i,\; \forall
i\in\{r-m,\ldots,r\}$. This finally proves $\mathrm{det}(\tau_1) \notin
\mathfrak{p}_i \; \forall i=1,\ldots,r$ after the sequence of base changes, i.e.
$\mathrm{det}(\tau_1)$ is then $R$-regular, that what was to be shown.\\
\vspace{0.5cm}
Secondly, I now want to prove:
\begin{lemm}
A codimension 2 Gorenstein algebra $A$ over a local CM ring $(R,\mathfrak{m},k)$ with
$2\notin \mathfrak{m}$ has a resolution
\[
\begin{CD}
0 \to R^n @>{-\beta^t\choose\alpha}>> R^{2n} @>(\alpha\:\beta)>> R^n \to A \to 0
\end{CD}
\]
which is also of Koszul module type, i.e.
$\mathrm{det}(\alpha),\:\mathrm{det}(\beta)$ is an $R-$regular sequence.
\end{lemm}
\begin{proof}
Taking into account the above remark that one can dispose of the assumption "$R$ a
domain" the fact that $A$ has a resolution with the symmetry property above is proven
in [Gra], thm. 3.3., so I have to show that $\exists$ a base change in $R^{2n}$ which
preserves the relation $\alpha \beta^t=\beta \alpha^t$ and in the new base
$\mathrm{det}(\alpha),\:\mathrm{det}(\beta)$ is a regular sequence. The punch line
to show this is as in the foregoing argument except that everything is a little
harder because one has to keep track of preserving the symmetry: Therefore let again
be  $\mathfrak{p}_1,\ldots,\mathfrak{p}_r$ the associated primes of $R$, and I show
that $\exists$ a sequence of $r$ base changes in $R^{2n}$ preserving the symmetry
and such that after the $m$th base change $(\ast)$ above holds, the case $m=0$ being
trivial. For the inductive step, suppose $\det(\alpha)\in\mathfrak{p}_{r-m}$ to
rule out a trivial case; I write
$[i_1,\ldots,i_{\nu};j_1,\ldots,j_{n-
\nu}]\equiv\det(\alpha_{i_1}\ldots\alpha_{i_{\nu}}\:\beta_{
j_1}\ldots\beta_{j_{n- \nu}})$. Call a minor $[i_1,\ldots,i_{\nu};j_1,\ldots,j_{n-
\nu}]$
\underline{good} iff $\{ i_1,\ldots,i_{\nu}\} \cap \{j_1,\ldots,j_{n- \nu}
\}=\emptyset$.\\ I want to find a good minor that does not belong to
$\mathfrak{p}_{r-m}$ (possibly after a base change in $R^{2n}$). Therefore suppose
all the good minors belong to
$\mathfrak{p}_{r-m}$. Since $\mathrm{grade}\:I_n((\alpha\:\beta))\ge 2$ by
Eisenbud-Buchsbaum acyclicity, $\exists$ a minor $\notin \mathfrak{p}_{r-m}$
(which is not good). For $n=1$ this is a contradiction since all minors are good, and
I can suppose $n>1$ in the process of finding a good minor. Now choose a minor
$[I_1,\ldots,I_k;J_1,\ldots,J_{n-k}]$ such that
\begin{itemize}
\item
$[I_1,\ldots,I_k;J_1,\ldots,J_{n-k}]\notin \mathfrak{p}_{r-m}$
\item
$\mathrm{card}( \{ I_1,\ldots,I_k \} \cap \{J_1,\ldots,J_{n-k} \} )=:M_0$  is minimal
among the minors which do not belong to $\mathfrak{p}_{r-m}$.
\end{itemize}
I want to perform a base change in $R^{2n}$ not destroying the symmetry such that in
the new base
$\exists$ a minor $[T_1,\ldots,T_{k-1};S_1,\ldots,S_{n-k+1}]$ such that
\begin{itemize}
\item
$[T_1,\ldots,T_{k-1};S_1,\ldots,S_{n-k+1}]\notin \mathfrak{p}_{r-m}$
\item
$\mathrm{card}( \{ T_1,\ldots,T_{k-1} \} \cap \{ S_1,\ldots,S_{n-k+1} \} )=M_0-1.$
\end{itemize}
Continuing this process $M_0$ steps (i.e. performing $M_0$ successive base changes) I
can find a good minor not contained in $\mathfrak{p}_{r-m}$.\\
Let now $[T_1,\ldots,T_{k-1};S_1,\ldots,S_{n-k+1}]$ be given. Choose $H \in \{
I_1,\ldots,I_k \} \cap \{ J_1,\ldots,J_{n-k} \}$ and $L \in \{ 1,\ldots,n \}
- \{ I_1,\ldots,I_k \} \cup \{ J_1,\ldots,J_{n-k} \}$ (both of which exist).
Now perform the base change in $R^{2n}$ which corresponds to adding $\alpha_H$ to
$\beta_L$ and $\alpha_L$ to $\beta_H$ (preserving the symmetry), and consider
\[
\det(\alpha_{I_1}\ldots\hat{\alpha_H}\ldots\alpha_{I_k}\;\beta_{J_1}\ldots\beta_H+\alpha_L
\ldots\beta_L+\alpha_H\ldots\beta_{J_{n-k}}),
\]
an $n\times n-$minor of the transformed matrix which I can write as\\
$[T_1,\ldots,T_{k-1};S_1,\ldots,S_{n-k+1}]$, where $\{ T_1,\ldots,T_{k-1} \} =\{
I_1,\ldots,I_k \}-\{ H\}$, $\{ S_1,\ldots,S_{n-k+1} \} = \{J_1,\ldots,J_{n-k} \}
\cup \{ L \}$ and obviously, $\mathrm{card}( \{ T_1,\ldots,T_{k-1} \}$\\ $\cap \{
S_1,\ldots,S_{n-k+1} \} )=M_0-1.$ I want to prove that this minor does not belong
to $\mathfrak{p}_{r-m}$. For this I show that in fact
\begin{eqnarray*}
[T_1,\ldots,T_{k-1};S_1,\ldots,S_{n-k+1}]=\pm
[I_1,\ldots,I_k;J_1,\ldots,J_{n-k}] \\+\mathrm{"residual\; terms"},
\end{eqnarray*}
where "residual terms"$\in \mathfrak{p}_{r-m}$. Using the additivity of the
determinant in each column I find that "residual terms" consists of 3 summands two of
which clearly belong to $\mathfrak{p}_{r-m}$ because
$[I_1,\ldots,I_k;J_1,\ldots,J_{n-k}]$ was chosen such that $\mathrm{card}( \{ I_1,\ldots,I_k \} \cap \{J_1,\ldots,J_{n-k} \}
)=:M_0$  was minimal among the minors of the matrix before the base change which did
not belong to
$\mathfrak{p}_{r-m}$, whereas the third summand is (up to sign)
\[
\det(\alpha_{I_1}\ldots\hat{\alpha_H}\ldots\alpha_{I_k}\:\alpha_L\:\beta_{J_1}\ldots
\hat{\beta_H}
\ldots\beta_{J_{n-k}}\:\beta_L).
\]
To show that the latter is in $\mathfrak{p}_{r-m}$ I apply the so-called
"Pl\"ucker relations":\vspace{0.3cm}

\itshape
Given an $M\times N-$matrix, $M\le 
 N,\; a_1,\ldots ,a_p,b_q,\ldots ,b_M,c_1,\ldots ,c_s\in \{ 1,\ldots ,N \} 
 ,\; s=M-p+q-1>M,\; t=M-p > 0 $, one has
\[
(P) \sum\limits_{ {i_1<\ldots <i_t \atop i_{t+1}<\ldots <i_s} \atop \{
i_1,\ldots ,i_s\} =\{ 1,\ldots ,s\} }
\sigma(i_1,\ldots,i_s)[a_1,\ldots ,a_p\: c_{i_1}\ldots c_{i_t}][c_{i_{t+1}}\ldots
c_{i_s}\: b_q\ldots b_M]=0
\]
where $\sigma(i_1,\ldots,i_s)$ is the sign of the permutation ${1,\ldots ,s \choose
i_1,\ldots ,i_s}$ (see e.g. [B-He], lemma 7.2.3, p. 308).
\vspace{0.3cm}

\upshape
In my situation, I let $M:=n,\;N:=2n,\;p:=n-2,\;q:=n+1,\;s:=n+1$ and for the
columns corresponding to the $a$'s above I choose the $n-2$ columns
\[
\alpha_{I_1},\: \alpha_{I_2},\ldots ,\hat{\alpha_H},\ldots ,\alpha_{I_k},\: 
\beta_{J_1}, \ldots ,\hat{\beta_H},\ldots ,\beta_{J_{n-k}} 
\]
(in this order), for the columns corresponding to the $b$'s I choose the empty set
(which is allowable here), and finally for the columns corresponding to the $c$'s the
$n+2$ columns
\begin{displaymath} 
\alpha_H,\: \beta_H,\: \alpha_L,\: \beta_L,\: \alpha_{I_1},\: \alpha_{I_2},\ldots
,\: \hat{\alpha_H},\ldots,\: \alpha_{I_k},\: \beta_{J_1},\ldots ,\:
\hat{\beta_H},\ldots ,\beta_{J_{n-k}} 
\end{displaymath}\\
Applying (P) one gets 6 nonvanishing summands, 4 of which (namely 
\begin{eqnarray*}\det(\alpha_{I_1}\ldots\hat{\alpha_H}\ldots\alpha_{I_k}\:\beta_{J_1}\ldots
\hat{\beta_H}
\ldots\beta_{J_{n-k}}\: \alpha_H\:\alpha_L)\cdot \mathrm{(a \; second\; factor)},\\
\det(\alpha_{I_1}\ldots\hat{\alpha_H}\ldots\alpha_{I_k}\:\beta_{J_1}\ldots
\hat{\beta_H}
\ldots\beta_{J_{n-k}}\: \alpha_H\:\beta_L)\cdot \mathrm{(a \; second\; factor)},\\
\det(\alpha_{I_1}\ldots\hat{\alpha_H}\ldots\alpha_{I_k}\:\beta_{J_1}\ldots
\hat{\beta_H}
\ldots\beta_{J_{n-k}}\: \beta_H\:\alpha_L)\cdot \mathrm{(a \; second\; factor)},\\
\det(\alpha_{I_1}\ldots\hat{\alpha_H}\ldots\alpha_{I_k}\:\beta_{J_1}\ldots
\hat{\beta_H}
\ldots\beta_{J_{n-k}}\: \beta_H\:\beta_L)\cdot \mathrm{(a \; second\;
factor)} ) \end{eqnarray*} 
 are in
$\mathfrak{p}_{r-m}$ by the defining minimality property of $[I_1,\ldots
,I_k;J_1,\ldots  ,J_{n-k}]$ above. The remaining 2 summands add up to (watch the
signs!)\\
$\pm
2 \det(\alpha_{I_1}\ldots\hat{\alpha_H}\ldots\alpha_{I_k}\:\alpha_L\:\beta_{J_1}\ldots
\hat{\beta_H}
\ldots\beta_{J_{n-k}}\:\beta_L) \cdot [I_1,\ldots ,I_k;J_1,\ldots  ,J_{n-k}]$
which therefore is also in $\mathfrak{p}_{r-m}$. But $[I_1,\ldots ,I_k;J_1,\ldots 
,J_{n-k}]\notin \mathfrak{p}_{r-m}$ and 2 is a unit in $R$, therefore
$\det(\alpha_{I_1}\ldots\hat{\alpha_H}\ldots\alpha_{I_k}\:\alpha_L\:\beta_{J_1}\ldots
\hat{\beta_H}
\ldots\beta_{J_{n-k}}\:\beta_L)\in\mathfrak{p}_{r-m}$ as desired, since
$\mathfrak{p}_{r-m}$ is prime.\\
Hence inductively, after $M_0$ base changes in $R^{2n}$, I can find a good minor of
the transformed matrix that is not in $\mathfrak{p}_{r-m}$. I assume $[1,\ldots,n]\in
\mathfrak{p}_{r-m}$. I can now define
\begin{eqnarray*}
l_1:=\mathrm{min}\{c :\:\exists s_1,\ldots,s_{n-1}\:\mathrm{with}\:
s_1<s_2<\ldots<s_{n-1}<c\\
\mathrm{and}\:[s_1,\ldots,s_{n-1},c]\notin\mathfrak{p}_{r-m} 
\mathrm{\underline{and}}\: [s_1,\ldots,s_{n-1},c]\: \mathrm{is \:
good} \}
\end{eqnarray*}
(then $n<l_1\le 2n$) and inductively,
\begin{eqnarray*}
l_i:=\mathrm{min}\{c :\:\exists s_1',\ldots,s_{n-i}'\:\mathrm{with}\:
s_1'<\ldots<s_{n-i}'<c<l_{i-1}<\ldots<l_1\\ 
\mathrm{and}\: [s_1',\ldots,s_{n-i}',c,l_{i-1},\ldots,l_1] \: \mathrm{is
\: good} \\
\mathrm{and}\:[s_1',\ldots,s_{n-i}',c,l_{i-1},\ldots,l_1]\notin\mathfrak{p}_{r-m}\}
.\end{eqnarray*} 
Then $[l_n,\ldots,l_1]\notin \mathfrak{p}_{r-m}$ which is good and can therefore be
written as
$[l_n,\ldots,l_1]=[l^{\alpha}_1,\ldots,l^{\alpha}_h;l^{\beta}_1,\ldots,l^{\beta}_{n-h}]$ 
with $\{ l^{\alpha}_1,\ldots,l^{\alpha}_h\} \cap \{
l^{\beta}_1,\ldots,l^{\beta}_{n-h} \} =\emptyset$. Choose
$b\in(\bigcap\limits_{i=r-m+1}^{r}\mathfrak{p}_i)\backslash\mathfrak{p}_{r-m}$ and
perform a base change in $R^{2n}$ (preserving the symmetry) by adding $b$ times the
$l_i^{\beta}$ column of $\beta$ to the $l_i^{\beta}$ column of $\alpha$, for
$i=1,\ldots,n-h$. Then
\[
[1,\ldots,n]_{\mathrm{new}}=[1,\ldots,n]_{\mathrm{old}}\pm
b^{n-h}[l_n,\ldots,l_1]_{\mathrm{old}}+b\mu,
\]
where $\mu\in\mathfrak{p}_{r-m}$ by the defining minimality property of the
$l$'s. Thus $[1,\ldots,n]_{\mathrm{new}}\notin\mathfrak{p}_i$ for
$i=r-m,\ldots,r$, which is the inductive step for the property ($\ast$).
Therefore after a sequence of base changes that preserve the symmetry
$\alpha\beta^t=\beta\alpha^t$, $\det(\alpha)$ can be made an $R-$regular
element.

Let's sum up: I have that $\det(\alpha)$ is a nonzerodivisor in $R$, and want
to prove that $\exists$ a base change in $R^{2n}$ preserving the symmetry and
leaving $\alpha$ unchanged (i.e. fixing the first $n$ basis vectors of
$R^{2n}$) such that in the new base $\det(\beta)$ is a nonzerodivisor in
$R/(\det(\alpha))$. The argument is almost identical to the preceding one. In
fact, let $\mathfrak{q}_1,\ldots ,\mathfrak{q}_s$ be the associated primes of
$R/(\det(\alpha))$ which are exactly the minimal prime ideals containing
$(\det(\alpha))$ because $R/(\det(\alpha))$ is CM ($R$ is CM and
$\det(\alpha)$ is $R-$regular). Then the part of the above proof starting with
"$\ldots$ the symmetry: Therefore let again
be  $\mathfrak{p}_1,\ldots,\mathfrak{p}_r$ the associated primes of $R$, and I show
that $\exists$ a sequence of $r$ base changes in $R^{2n}$ $\ldots$" and ending
with "$\ldots$ Choose $H \in \{
I_1,\ldots,I_k \} \cap \{ J_1,\ldots,J_{n-k} \}$ and $L \in \{ 1,\ldots,n \}
- \{ I_1,\ldots,I_k \} \cup \{ J_1,\ldots,J_{n-k} \}$ $\ldots$" goes through
\underline{verbatim} (and has to be inserted here) if throughout one replaces
$r$ with $s$, $\det(\alpha)$ with $\det(\beta)$, and the symbol
"$\mathfrak{p}$" with "$\mathfrak{q}$". Thereafter, a slight change is
necessary because in the process of finding a good minor, i.e. in the course
of the $M_0$ base changes on $R^{2n}$ that transform $(\alpha\:\beta)$ s.t. in
the new base $\exists$ a good minor, the shape of $\beta$ is changed. This
change must preserve the property $\det(\beta)\notin\mathfrak{q}_1,\ldots
,\mathfrak{q}_{s-m+1}$ in order not to destroy the induction hypothesis. The
way out is as follows:\\
Choose $\zeta
\in (\bigcap\limits_{i=r-m+1}^{r}\mathfrak{q}_i)\backslash\mathfrak{q}_{r-m}$,
which is possible since the $q$'s all have height 1. Now perform the base change in $R^{2n}$ which corresponds to adding
$\zeta\alpha_H$ to
$\beta_L$ and $\zeta\alpha_L$ to $\beta_H$ (preserving the symmetry), and consider
\[
\det(\alpha_{I_1}\ldots\hat{\alpha_H}\ldots\alpha_{I_k}\;\beta_{J_1}\ldots\beta_H+\zeta\alpha_L
\ldots\beta_L+\zeta\alpha_H\ldots\beta_{J_{n-k}}),
\]
an $n\times n-$minor of the transformed matrix which I can write as\\
$[T_1,\ldots,T_{k-1};S_1,\ldots,S_{n-k+1}]$, where $\{ T_1,\ldots,T_{k-1} \} =\{
I_1,\ldots,I_k \}-\{ H\}$, $\{ S_1,\ldots,S_{n-k+1} \} = \{J_1,\ldots,J_{n-k} \}
\cup \{ L \}$ and obviously, $\mathrm{card}( \{ T_1,\ldots,T_{k-1} \}$\\ $\cap \{
S_1,\ldots,S_{n-k+1} \} )=M_0-1.$ I want to prove that this minor does not belong
to $\mathfrak{q}_{s-m}$ and furthermore that 
\[
\det(\beta_1\ldots \beta_H+\zeta\alpha_L \ldots \beta_L+\zeta\alpha_H \ldots \beta_n) \notin\mathfrak{q}_1,\ldots
,\mathfrak{q}_{s-m+1}.
\] 
The latter statement is obvious by the choice of $\zeta$ (and multilinearity of
determinants). The former one follows if I show
\begin{eqnarray*}
[T_1,\ldots,T_{k-1};S_1,\ldots,S_{n-k+1}]=\pm
\zeta [I_1,\ldots,I_k;J_1,\ldots,J_{n-k}]\\ +\mathrm{"residual\; terms"},
\end{eqnarray*}
where "residual terms"$\in \mathfrak{q}_{s-m}$ because $\zeta$ and
$[I_1,\ldots ,I_k;J_1,\ldots ,J_{n-k}]$ are both $\notin\mathfrak{q}_{s-m}$ by
assumption. Again "residual terms" consists of 3 summands two of which belong
to $\mathfrak{q}_{s-m}$ because of the defining minimality property of
$[I_1,\ldots,I_k;J_1,\ldots,J_{n-k}]$. The third summand is up to sign 
\[
\zeta \det(\alpha_{I_1}\ldots\hat{\alpha_H}\ldots\alpha_{I_k}\:\alpha_L\:\beta_{J_1}\ldots
\hat{\beta_H}
\ldots\beta_{J_{n-k}}\:\beta_L),
\] 
therefore it suffices to show
$\det(\alpha_{I_1}\ldots\hat{\alpha_H}\ldots\alpha_{I_k}\:\alpha_L\:\beta_{J_1}\ldots
\hat{\beta_H}
\ldots\beta_{J_{n-k}}\:\beta_L)\in \mathfrak{q}_{s-m}$. This is done word by word
as in the passage of the first part of this proof starting with "$\ldots$ I apply
the so-called "Pl\"ucker relations":$\ldots$ " and ending "$\ldots$ since
$\mathfrak{p}_{r-m}$ is prime.$\ldots$ ", taking into account the afore-mentioned
changes in notation.\\
The rest of the proof is as follows: Inductively, I can find a good minor of the
transformed matrix that is not in $\mathfrak{q}_{s-m}$. To avoid a trivial case,
I assume $[n+1,\ldots ,2n]\in \mathfrak{q}_{s-m}$. Now I define
\begin{eqnarray*}
L_1:=\mathrm{max}\{c :\:\exists s_2,\ldots,s_n\:\mathrm{with}\:
c<s_2<s_3<\ldots<s_n\\ \mathrm{and}\:[c,s_2,\ldots,s_n]\notin\mathfrak{q}_{s-m}
 \mathrm{\underline{and}}\: [c,s_2,\ldots,s_n]\: \mathrm{is \:
good} \}
\end{eqnarray*}
(then $1\le l_1< n+1$) and inductively,
\begin{eqnarray*}
L_i:=\mathrm{max}\{c :\:\exists s_{i+1}',\ldots,s_{n}'\:\mathrm{with}\:
L_1<\ldots<L_{i-1}<c<s_{i+1}'<\ldots<s_n'\\
\mathrm{and}\: [L_1,\ldots,L_{i-1},c,s_{i+1}',\ldots,s_n'] \: \mathrm{is
\: good} \\
\mathrm{and}\:[L_1,\ldots,L_{i-1},c,s_{i+1}',\ldots,s_n']\notin\mathfrak{q}_{s-m}\}
.\end{eqnarray*} 
Then  $[L_1,\ldots,L_n]\notin \mathfrak{q}_{s-m}$ and is good (furthermore
$L_n >n$ since $[1,\ldots,n] \in \mathfrak{q}_{s-m}$). I can write
$[L_1,\ldots,L_1]=[L^{\alpha}_1,\ldots,L^{\alpha}_h;L^{\beta}_1,\ldots,L^{\beta}_{n-h}]$ 
with\\ $\{ L^{\alpha}_1,\ldots,L^{\alpha}_h\} \cap \{
L^{\beta}_1,\ldots,L^{\beta}_{n-h} \} =\emptyset$. Choose
$b\in(\bigcap\limits_{i=s-m+1}^{s}\mathfrak{q}_i)\backslash\mathfrak{q}_{s-m}$ and
perform a base change in $R^{2n}$ (preserving the symmetry) by adding $b$ times the
$L_i^{\alpha}$ column of $\alpha$ to the $L_i^{\alpha}$ column of $\beta$, for
$i=1,\ldots,h$. Then
\[
[n+1,\ldots,2n]_{\mathrm{new}}=[n+1,\ldots,2n]_{\mathrm{old}}\pm
b^{h}[L_1,\ldots,L_n]_{\mathrm{old}}+b\mu,
\]
where $\mu\in\mathfrak{q}_{s-m}$ by the defining maximality property of the
$L$'s. Thus $[n+1,\ldots,2n]_{\mathrm{new}}\notin\mathfrak{q}_i$ for
$i=s-m,\ldots,s$, which is the inductive step. Therefore
after a sequence of base changes that preserve the symmetry
$\alpha\beta^t=\beta\alpha^t$ (and leave $\det(\alpha)$ unaltered) $\det(\beta)$ can
be made an
$R/(\det(\alpha))-$regular element, i.e. $\det(\alpha),\:\det(\beta)$ is an
$R-$regular sequence, which proves the lemma.      
\end{proof}
\vspace{3cm}
\begin{center} \bfseries References\end{center}
\mdseries
\small
\begin{itemize}
\item[\textrm{[A-C-G-H]}] E. Arbarello, M. Cornalba, P. Griffiths, J. Harris:
\emph{Geometry of Algebraic Curves}, Springer-Verlag, New York, N.Y. (1985) 

\item[\textrm{[Beau]}] A. Beauville: \emph{Complex algebraic surfaces}, London
Mathematical Society Student Texts 34, CUP (1996)

\item[\textrm{[Bom]}] E. Bombieri: \emph{Canonical models of surfaces of general type},
I.H.E.S. Publ. Math. 42 (1973), p. 171-219

\item[\textrm{[B-He]}] W. Bruns, J. Herzog: \emph{Cohen-Macaulay rings}, CUP (1998)

\item[\textrm{[B-V]}] W. Bruns, U. Vetter: \emph{Determinantal rings}, Springer LNM
1327 (1988)

\item[\textrm{[Cat1]}] F. Catanese: \emph{Babbage's conjecture, contact of
surfaces, symmetric determinantal varieties and applications}, Inv.
Math. 63 (1981), pp.433-465

\item[\textrm{[Cat1b]}] F. Catanese: \emph{On the moduli spaces of surfaces of general
type}, J. Differential Geom., 19 (1984), p. 483-515
 
\item[\textrm{[Cat2]}] F. Catanese: \emph{Commutative
algebra methods and equations of regular surfaces}, Algebraic
Geometry-Bucharest 1982, Springer LNM 1056 (1984), pp.68-111

\item[\textrm{[Cat3]}] F.Catanese: \emph{Equations of pluriregular varieties of general
type}, Geometry today-Roma 1984, Progr. in Math. 60, Birkh\"auser
(1985), pp. 47-67)

\item[\textrm{[Cat4]}] F. Catanese: \emph{Homological Algebra and Algebraic
surfaces}, Proc. Symp. in Pure Math., Volume 62.1, 1997, pp. 3-56

\item[\textrm{[Cil]}] C. Ciliberto: \emph{Canonical
surfaces with
$p_g=p_a=5$ and
$K^2=10$}, Ann. Sc.  Norm. Sup. s.IV, v.IX, 2 (1982), pp.287-336

\item[\textrm{[C-E-P]}] C. de Concini, D. Eisenbud, C. Procesi: \emph{Hodge algebras},
Ast\'{e}risque 91 (1982)

\item[\textrm{[C-S]}] C. de Concini, E. Strickland: \emph{On the variety of complexes},
Adv. in Math. 41 (1981), pp. 57-77 

\item[\textrm{[D-E-S]}] W. Decker, L. Ein, F.-O. Schreyer: \emph{Construction of
surfaces in $\mathbb{P}^4$}, J. Alg. Geom. 2 (1993), pp.185-237

\item[\textrm{[Ei]}] D. Eisenbud:  \emph{Commutative Algebra with a view towards
Algebraic Geometry}, Springer G.T.M. 150, New York (1995)

\item[\textrm{[E-U]}] D. Eisenbud, B. Ulrich: \emph{Modules that are Finite
Birational Algebras}, Ill. Jour. Math. 41, 1 (1997), pp. 10-15

\item[\textrm{[En]}] F. Enriques: \emph{Le superficie algebriche},
Zanichelli, Bologna, 1949

\item[\textrm{[Gra]}] M. Grassi:
\emph{Koszul modules and Gorenstein algebras}, J. Alg. 180 (1996),
pp.918-953

\item[\textrm{[Gri-Ha]}] P. Griffiths, J. Harris: \emph{Principles of Algebraic
Geometry}, Wiley, New York (1978)

\item[\textrm{[Ha]}] J. Harris: \emph{Algebraic Geometry, A First Course}, Springer
G.T.M. 133, New York (1993)

\item[\textrm{[Har]}] R. Hartshorne: \emph{Algebraic Geometry}, Springer G.T.M. 52, New
York (1977)

\item[\textrm{[dJ-vS]}] T. de Jong, D. van Straten: \emph{Deformation of
the normalization of hypersurfaces}, Math. Ann. 288 (1990), S.527-547

\item[\textrm{[Lip]}] J. Lipman: \emph{Dualizing sheaves, differentials and residues on
algebraic varieties}, Ast\'{e}risques 117 (1984)

\item[\textrm{[M-P]}] D. Mond, R. Pellikaan: \emph{Fitting
ideals and multiple points of analytic mappings}, Springer LNM 1414
(1987), pp.107-161

\item[\textrm{[P-S1]}] C. Peskine, L. Szpiro: \emph{Dimension projective finie et
cohomologie locale}, Publ. Math. IHES 42 (1972)

\item[\textrm{[P-S2]}] C. Peskine, L. Szpiro: \emph{Liaison des vari\'{e}t\'{e}s
alg\'{e}briques}, Inv. Math. 26 (1972), pp. 271-302

\item[\textrm{[Ro\ss ]}] D. Ro\ss berg: \emph{Kanonische Fl\"achen mit $p_g=5,\:q=0$
und
$K^2=11,12$}, Doktorarbeit Bayreuth, November 1996

\item[\textrm{[Sern]}] E. Sernesi: \emph{L'unirazionalit\`a della variet\`a dei 
moduli delle curve di genere dodici}, Ann. Scuola Norm. Pisa 8 (1981), pp.
405-439

\item[\textrm{[Wal]}] C. Walter: \emph{Symmetric sheaves and Witt groups}, preprint

\end{itemize}
    
\end{document}